\input amstex\documentstyle {amsppt}  
\pagewidth{12.5 cm}\pageheight{19 cm}\magnification\magstep1
\topmatter
\title Character sheaves on disconnected groups, V\endtitle
\author G. Lusztig\endauthor
\address Department of Mathematics, M.I.T., Cambridge, MA 02139\endaddress
\thanks Supported in part by the National Science Foundation\endthanks
\endtopmatter   
\document
\define\mpb{\medpagebreak}
\define\co{\eta}

\define\bSi{\bar\Si}
\define\bH{\bar H}

\define\bboc{\bar\boc}

\define\hZ{\hat Z}
\define\hK{\hat K}

\define\uce{\un{\ce}}

\define\po{\text{\rm pos}}

\define\si{\sim}
\define\wt{\widetilde}
\define\sqc{\sqcup}

\define\hA{\hat A}

\define\hY{\hat Y}

\define\bj{\bar j}

\define\bY{\bar Y}

\define\lb{\linebreak}

\define\op{\oplus}

\define\em{\emptyset}
\define\imp{\implies}

\define\n{\notin}
\define\iy{\infty}
\define\m{\mapsto}
\define\do{\dots}

\define\bsl{\backslash}

\define\lra{\leftrightarrow}

\define\sub{\subset}
\define\bxt{\boxtimes}
\define\T{\times}
\define\ti{\tilde}
\define\nl{\newline}
\redefine\i{^{-1}}

\define\un{\underline}

\define\ot{\otimes}
\define\bbq{\bar{\QQ}_l}

\define\Ad{\text{\rm Ad}}
\define\Hom{\text{\rm Hom}}
\define\End{\text{\rm End}}

\define\ind{\text{\rm ind}}

\define\res{\text{\rm res}}

\define\tr{\text{\rm tr}}

\define\supp{\text{\rm supp}}

\define\a{\alpha}
\redefine\b{\beta}
\redefine\c{\chi}
\define\g{\gamma}
\redefine\d{\delta}
\define\e{\epsilon}
\define\et{\eta}
\define\io{\iota}
\redefine\o{\omega}
\define\p{\pi}
\define\ph{\phi}
\define\ps{\psi}
\define\r{\rho}
\define\s{\sigma}
\redefine\t{\tau}
\define\th{\theta}

\define\z{\zeta}
\define\x{\xi}

\define\vt{\vartheta}

\redefine\G{\Gamma}
\redefine\D{\Delta}

\define\Si{\Sigma}
\define\Th{\Theta}

\define\Ph{\Phi}

\redefine\aa{\bold a}

\define\boc{\bold c}
\define\dd{\bold d}

\define\kk{\bold k}

\define\nn{\bold n}

\redefine\ss{\bold s}
\redefine\tt{\bold t}

\redefine\AA{\bold A}

\define\EE{\bold E}
\define\FF{\bold F}

\define\II{\bold I}

\define\NN{\bold N}

\define\QQ{\bold Q}

\define\UU{\bold U}

\define\WW{\bold W}
\define\ZZ{\bold Z}

\define\ca{\Cal A}
\define\cb{\Cal B}

\define\cd{\Cal D}
\define\ce{\Cal E}
\define\cf{\Cal F}

\define\ch{\Cal H}

\define\cl{\Cal L}
\define\cm{\Cal M}

\define\cp{\Cal P}

\define\cs{\Cal S}
\define\ct{\Cal T}
\define\cu{\Cal U}
\define\cv{\Cal V}
\define\cw{\Cal W}
\define\cz{\Cal Z}
\define\cx{\Cal X}

\define\fA{\frak A}

\define\fC{\frak C}
\define\fD{\frak D}

\define\fK{\frak K}

\define\fT{\frak T}

\define\fZ{\frak Z}
\define\fX{\frak X}

\define\ta{\ti a}

\define\tf{\ti f}

\define\tj{\ti j}

\define\tp{\ti p}

\define\tA{\ti A}

\define\tG{\ti G}

\define\tK{\ti K}

\define\tM{\ti M}

\define\tP{\ti P}  
\define\tQ{\ti Q}

\define\tS{\ti S}

\define\tY{\ti Y}

\define\sh{\sharp}

\define\bS{\bar S}

\define\bce{\bar\ce}
\define\tce{\ti\ce}
\define\bul{\bullet}

\define\bae{\bar\e}
\define\bacv{\bar\cv}
\define\hce{\hat{\ce}}
\define\che{\check}
\define\BE{B}
\define\BBD{BBD}
\define\CS{L3}
\define\AD{L9}
\define\PCS{L10}
\head Introduction\endhead
Throughout this paper, $G$ denotes a fixed, not necessarily connected, reductive
algebraic group over an algebraically closed field $\kk$. This paper is a part of a
series \cite{\AD} which attempts to develop a theory of character sheaves on $G$. The 
numbering of the sections and references continues that of the earlier Parts.

Section 23 is a generalization of results in \cite{\CS, II,\S7}. It is a preparation
for the proof of the orthogonality formulas for certain characteristic functions
in Section 24 which generalize those in \cite{\CS, II,\S9,\S10}. Section 25 
describes the cohomology sheaves of a class of complexes which includes the admissible 
complexes on $G$. In particular we show that these cohomology sheaves restricted to
any stratum of $G$ are local systems of a particular kind. In the connected case this
reduces to a strengthening of \cite{\CS, III,14.2(a)}. In Section 26 we give a variant 
of the definition of parabolic character sheaves in \cite{\PCS} in terms of admissibile
complexes. Note that even if one is only interested in parabolic character sheaves of
connected groups, one cannot avoid using the theory of character sheaves on disconnected 
groups. In Section 27 we discuss the induction functor. The present treatment
differs from one in the connected case, given in \cite{\CS, I,\S4}.

{\it Notation.} Let $\cf$ be a local system on an algebraic variety $Y$. If $\kk$ is an
algebraic closure of a finite field $\FF_q$, $F:Y@>>>Y$ is the Frobenius map for an 
$\FF_q$-rational structure and $\e:F^*\cf@>\si>>\cf$ is an isomorphism, we denote by 
$\che\e:F^*\che\cf@>>>\che\cf$ the unique isomorphism such that for any $y\in Y$, 
$\che\e:\che\ce_{F(y)}@>>>\che\ce_y$ is the isomorphism transpose inverse to 
$\e:\ce_{F(y)}@>>>\ce_y$.

If $X$ is an algebraic variety and $A\in\cd(X)$ is $IC(X',\cf)$ extended by $0$ on 
$X-X'$ where $X'$ is a closed irreducible subvariety of $X$ and $\cf$ is a local system
on an open dense smooth subvariety $X'_0$ of $X'$, let $\che A=IC(X',\che\cf)$ extended
by $0$ on $X-X'$. We have also $\che A=\fD(A)[-2\dim X']$. If $\kk$ is an algebraic 
closure of a finite field $\FF_q$, $F:X@>>>X$ is the Frobenius map for an 
$\FF_q$-rational structure such that $F(X')=X',F(X'_0)=X'_0$ and $\a:F^*A@>\si>>A$ is 
an isomorphism which restricts to $\e:F^*\cf@>\si>>\cf$ over $X_0$, we denote by 
$\che\a:F^*\che A@>>>\che A$ the unique isomorphism which restricts to 
$\che\e:F^*\che\cf@>\si>>\che\cf$ over $X_0$.

If $f:X@>>>Y$ is a smooth morphism between algebraic varieties with connected fibres of
dimension $\d$, we set $\tf A=f^*A[\d]$ for any $A\in\cd(Y)$.
Let $\cd(X)^{\le 0}$ be as in \cite{\CS, I,1.3}.
Let $\cm(X)$ be the category of perverse sheaves on $X$.

If $D$ is a connected component of $G$, a simple perverse sheaf $A$ on $D$ is said to 
be admissible if $A$, regarded as a simple perverse sheaf on $G$, zero on $G-D$, is 
admissible in the sense of 6.7.

\head Contents\endhead
23. Strongly cuspidal local systems.

24. Orthogonality.

25. Properties of cohomology sheaves.

26. The variety $Z_{J,D}$.

27. Induction.

\head 23. Strongly cuspidal local systems\endhead
\subhead 23.1\endsubhead
Let $S$ be an isolated stratum of $G$ and let $\ce\in\cs(S)$. Let $D$ be the connected 
component of $G$ that contains $S$.

If $P$ is a parabolic of $G^0$ such that $S\sub N_GP$ and $R$ is an $U_P$-coset in 
$N_GP$, let $d_R$ be $\dim S-\dim{}^D\cz_{G^0}^0$ minus the dimension of the $P/U_P$ 
orbit of $R/U_P$ in $N_GP/P$. 

We show that the following two conditions are equivalent:

(i) {\it $\ce$ is a cuspidal local system;}

(ii) {\it for any $P,R$ as above such that $P\ne G^0$ and any irreducible component 
$\fC$ of $S\cap R$ of dimension $d_R/2$, the restriction of $\che\ce$ to some/any 
smooth open dense subset of $\fC$ has no direct summand $\bbq$.}
\nl
Let $\fC_0$ be a smooth open dense subset of $\fC$. We have 
$H^{d_R}_c(\fC,\ce)\cong H^{d_R}_c(\fC_0,\ce)$ and, by Poincar\'e duality, the last 
vector space is isomorphic to $H^0(\fC_0,\che\ce)$, a vector space whose dimension is
the multiplicity of $\bbq$ in a decomposition of $\che\ce|_{\fC_0}$ as a direct sum of
irreducible local systems. It remains to note that, by 6.2, 
$H^{d_R}_c(S\cap R,\ce)\cong\op_\fC H^{d_R}_c(\fC,\ce)$ where $\fC$ runs over the 
irreducible components of $S\cap R$ of dimension $d_R/2$.

Since condition (ii) for $\ce$ is equivalent to condition (ii) for $\che\ce$, we 
deduce:

(a) {\it $\ce$ is cuspidal if and only if $\che\ce$ is cuspidal.}

\subhead 23.2\endsubhead
Let $S,\ce,D,P$ be as in 23.1; assume that $P\ne G^0$. Let $K_0=IC(\bS,\ce)$. We 
show:

(a) {\it if $R,d_R$ are as in 23.1 and $H^i_c((\bS-S)\cap R,K_0)\ne 0$ then $i<d_R$.}
\nl
Indeed, there exist $i',i''$ and a stratum $S_1\sub\bS-S$ with $i=i'+i''$ and 
$H^{i'}_c(S_1\cap R,\ch^{i''}K_0)\ne 0$. By 6.2 we have 
$\dim(S_1\cap R)\le(\dim S_1-\dim S+d_R)/2$. Since $i'\le 2\dim(S_1\cap R)$ we have 
$i'\le\dim S_1-\dim S+d_R$. Since $\ch^{i''}K_0\ne 0$ at all points of $S_1$, we have
$i''<\dim S-\dim S_1$. Hence $i=i'+i''<\dim S_1-\dim S+d_R+\dim S-\dim S_1$ and 
$i<d_R$, as required.

Next we show that:

(b) {\it if $R,d_R$ are as in 23.1, the natural map 
$H^i_c(S\cap R,\ce)@>j_i>>H^i_c(\bS\cap R,K_0)$ is an isomorphism for $i\ge d_R$.}
\nl
Indeed, from (a) we see that $j_i$ is surjective for $i=d_R$ and an isomorphism for 
$i>d_R$. Moreover, for $i=d_R$, the kernel of $j_i$ equals the image of the natural map
$f:H^{d_R-1}_c((\bS-S)\cap R,K_0)@>>>H^{d_R}_c(S\cap R,\ce)$. It is enough to show that
$f=0$. We argue as in the proof of 8.3(b). We may assume that $\kk$ is an algebraic 
closure of a finite field $\FF_q$, that $G$ has a fixed $\FF_q$-structure with 
Frobenius map $F:G@>>>G$, that $P,R$ and any stratum in $\bS$ are defined over $\FF_q$ 
and that we are given an isomorphism $F^*\ce@>\si>>\ce$ which makes $\ce$ into a local 
system of pure weight $0$. Then $H^{d_R-1}_c((\bS-S)\cap R,K_0),H^{d_R}_c(S\cap R,\ce)$
have natural (Frobenius) endomorphisms compatible with $f$. To show that $f=0$ it is 
enough to show that

$(*)$ $H^{d_R}_c(S\cap R,\ce)$ is pure of weight $d_R$;

$(**)$ $H^{d_R-1}_c((\bS-S)\cap R,K_0)$ is mixed of weight $\le d_R-1$.
\nl
Now $(*)$ is clear since $\dim(S\cap R)\le d_R/2$ (see 6.2). We prove $(**)$. As in the
proof of (a), it is enough to show that for any stratum $S_1\sub\bS-S$ and any $i',i''$
such that $i'+i''=d_R-1$, $H^{i'}_c(S_1\cap R,\ch^{i''}K_0)$ is mixed of weight 
$\le d_R-1$. By Gabber's theorem \cite{\BBD, 5.3.2}, the local system $\ch^{i''}K_0$ on
$S_1$ is mixed of weight $\le i''$. Using Deligne's theorem \cite{\BBD, 5.1.14(i)} we 
deduce that $H^{i'}_c(S_1\cap R,\ch^{i''}K_0)$ is mixed of weight $\le i'+i''=d_R-1$. 
This completes the proof of (b).

\mpb

We show that the following two conditions for $S,\ce,P$ are equivalent:

(i) {\it for any $U_P$-coset $R$ in $N_GP$ we have $H^{d_R}_c(S\cap R,\ce)=0$,
$d_R$ as in 23.1;}

(ii) {\it for any $i$, the set $\cx_i$ consisting of the $U_P$-cosets $R$ in $N_GP$ 
such that $H^i_c(\bS\cap R,K_0))\ne 0$ has dimension $<\dim S-i$.}
\nl
Consider the set of all $U_P$-cosets $R$ in $N_GP$ such that $\bS\cap R\ne\em$. On this
set we have a ${}^D\cz_{G^0}^0\T P/U_P$-action $(z,p):R\m zpRp\i$; this action has only
finitely many orbits and some orbit has dimension $\dim S-d_R$. Also, this action leaves stable each of the subsets $\cx_i$ in (ii).

Assume that (ii) holds. Let $R$ be an $U_P$-coset in $N_GP$. Assume that \lb
$H^{d_R}_c(S\cap R,\ce)\ne 0$. By (b) we have $R\in\cx_i$ with $i=d_R$. Hence the orbit
of $R$ is contained in $\cx_i$. It follows that $\dim\cx_i\ge\dim S-d_R$ contradicting 
(ii). Thus, (i) holds.

Conversely, assume that (i) holds. To establish (ii), we assume that $i$ is such that
$\cx_i\ne\em$; it is enough to show that for any ${}^D\cz_{G^0}^0\T P/U_P$-orbit $\co$ 
in $\cx_i$ we have $\dim\co<\dim S-i$. If $R\in\co$ we have $\dim\co=\dim S-d_R$. Hence
it is enough to show that $\dim S-d_R<\dim S-i$ or that $i<d_R$. Assume that 
$i\ge d_R$. Since $\cx_i\ne\em$ we have $H^i_c(\bS\cap R,K_0)\ne 0$ for some $R$ hence
by (b) we have $H^i_c(S\cap R,\ce)\ne 0$. If $i>d_R$ this is impossible, by 6.2. If 
$i=d_R$ this is impossible since (i) holds. Thus we have $i<d_R$ and (ii) holds.

\subhead 23.3\endsubhead
In the setup of 23.1 we say that $\ce$ is {\it clean} if $IC(\bS,\ce)|_{\bS-S}=0$. We 
say that $\ce\in\cs(S)$ is {\it strongly cuspidal} if for any parabolic $P$ of $G^0$ 
such that $P\ne G^0,S\sub N_GP$ and any $U_P$-coset $R$ in $N_GP$ we have
$H^i_c(\bS\cap R,K_0)=0$ for all $i$.

From the equivalence of (i),(ii) in 23.2 we see that:

{\it if $\ce$ is strongly cuspidal then $\ce$ is cuspidal.}
\nl
If $\ce$ is assumed to be clean, then the condition that $\ce$ is strongly cuspidal is
equivalent to the following condition: for any parabolic $P$ of $G^0$ such that 
$P\ne G^0,S\sub N_GP$ and any $U_P$-coset $R$ in $N_GP$ we have $H^i_c(S\cap R,\ce)=0$
for all $i$.

\mpb

Let $D$ be a connected component of $G$ and let $P$ be a parabolic of $G^0$ such that 
that $N_GP\cap D\ne\em$. Let $L$ be the Levi of $P$. Let $G'=N_GP\cap N_GL$, a 
reductive group with $G'{}^0=L$. Let $D'=G'\cap D$, a connected component of $G'$. 
Define a homomorphism $\p:N_GP@>>>G'$ by $\p(z\o)=z$ where $z\in N_GP\cap N_GL$,
$\o\in U_P$ (see 1.26). The restriction of $\p$ to $N_GP\cap D@>>>D'$ is denoted again
by $\p$. Let $i:N_GP\cap D@>>>D$ be the inclusion. We define a functor 
$$\res_D^{D'}:\cd(D)@>>>\cd(D')$$ 
by $\res_D^{D'}A=\p_!i^*A$.

Let $A$ be a perverse sheaf on $D$. We say that $A$ is {\it cuspidal} if
$$\res_D^{D'}A[-1]\in\cd(D')^{\le 0}$$
\nl
for any $P,L,D'$ as above with $P\ne G^0$. We say that $A$ is {\it strongly cuspidal} 
if 
$$\res_D^{D'}(A)=0\in\cd(D')$$ 
for any $P,L,D'$ as above with $P\ne G^0$. Clearly, if $A$ is strongly cuspidal, then 
it is cuspidal.

Assume now that $S$ is an isolated stratum in $D$ and $\ce\in\cs(S)$. Let 
$A=IC(\bS,\ce)[\dim S]$ regarded as a perverse sheaf on $D$, zero on $D-\bS$. Clearly, 
$A$ is strongly cuspidal if and only if $\ce$ is a strongly cuspidal local system. We 
show:

(a) {\it $A$ is cuspidal if and only if $\ce$ is a cuspidal local system.}
\nl
By the equivalence of (i),(ii) in 23.2, the condition that $\ce$ is a cuspidal local
system is that, for any $P,L,D'$ be as above with $P\ne G^0$ and any $j\in\ZZ$ we have
$$\dim(\supp\ch^j(\res_D^{D'}A[-\dim S]))<\dim S-j$$
or equivalently $\dim(\supp\ch^j(\res_D^{D'}A))<-j$ that is,
$\dim(\supp\ch^j(\res_D^{D'}A[-1]))\le-j$. This is the same as
$\res_D^{D'}A[-1]\in\cd(D')^{\le 0}$. This proves (a).

\subhead 23.4\endsubhead
Let $\ce\in\cs(S)$. Let $s\in S_s$ and let $G'=Z_G(s)$. Let $\boc$ be a 
$G'{}^0$-conjugacy class in $\cv=\{v\in G';v\text{ unipotent, }sv\in S\}$. Let $\d$ be 
the connected component of $G'$ that contains $\boc$. Let $\ce'$ be the inverse image 
of $\ce$ under the map ${}^\d\cz_{G'{}^0}^0\boc@>>>S, g\m sg$. Then $\ce'\in\cs(S')$ 
where $S'={}^\d\cz_{G'{}^0}^0\boc$. Let $K_0=IC(\bS,\ce),K'_0=IC(\bS',\ce')$. We show:

(a) {\it $\ce$ is clean (with respect to $G$) if and only if $\ce'$ is clean (with 
respect to $G'$).}
\nl
Let $\ce''$ be the inverse image of $\ce$ under $\cv@>>>S,v\m sv$. Since $\cv$ is a
smooth equidimensional variety (it admits a transitive action of an algebraic group), 
the complex $IC(\bar\cv,\ce'')\in\cd(\bar\cv)$ is well defined.

Using 1.22 we see that $\p:\bS@>>>S_s,g\m g_s$ is a morphism of algebraic varieties.
Now ${}^D\cz_{G^0}^0\T G^0$ acts on $\bS$ and $S_c$ compatibly with $\p$ so that the 
action on $S_c$ is transitive. Since the fibre $\p\i(s)$ may be identified with 
$\bacv$, we see that $IC(\bacv,\ce'')=h^*K_0$ where $h:\bacv@>>>\bS,v\m sv$. It follows
that $\ce$ is clean if and only if $IC(\bacv,\ce'')|_{\bacv-\cv}=0$. We have 
$\cv=\sqc_{i\in[1,m]}\boc_i$ where $\boc_i$ are $G'{}^0$-conjugacy classes. Hence 
$IC(\bacv,\ce'')=\op_{i\in[1,m]}IC(\bboc_i,\ce''|_{\bboc_i})$. Thus, 
$IC(\bacv,\ce'')|_{\bacv-\cv}=0$ if and only if
$IC(\bboc_i,\ce''|_{\bboc_i})|_{\bboc_i-\boc_i}=0$ for all $i$. By the homogeneity of 
$\cv$ this is equivalent to the condition that 
$IC(\bboc,\ce''|_{\bboc})|_{\bboc-\boc}=0$. (We have $\boc=\boc_i$ for some $i$.) This 
last condition is clearly equivalent to the condition that $K'_0|_{\bS'-\bS'}=0$. This
proves (a).

We show:

(b) {\it If $\ce$ is strongly cuspidal (with respect to $G$) then $\ce'$ is strongly 
cuspidal (with respect to $G'$).}
\nl
Let $Q$ be a parabolic of $G'{}^0$ such that $Q\ne G'{}^0,S'\sub N_{G'}Q$. We must show
that for any $z\in{}^\d\cz_{G'{}^0}^0,u\in\bboc\cap Q$ we have $H^i_c(zuU_Q,K'_0)=0$ 
for all $i$. We may assume that $z=1$. Hence we must show that 
$H^i_c(uU_Q\cap\bboc,K'_0)=0$ for all $i$. Since $H^i_c(uU_Q\cap\bboc,K'_0)$ is a 
direct summand of $H^i_c(uU_Q\cap\bacv,IC(\bacv,\ce''))$ it is enough to show that
$H^i_c(uU_Q\cap\bacv,IC(\bacv,\ce''))=0$ for all $i$. Since $IC(\bacv,\ce'')=h^*K_0$,  
it is enough to show that $H^i_c(suU_Q\cap s\bacv,K_0)=0$ for all $i$.
We have $suU_Q\cap s\bacv=suU_Q\cap\bS$ hence it is enough to show that 
$H^i_c(suU_Q\cap\bS,K_0)=0$ for all $i$. By 1.18(a) we can find a parabolic $P$ of 
$G^0$ such that $P\cap G'{}^0=Q$ and $su\in N_GP$. Clearly, $P\ne G^0$. Let 
$f:suU_P\cap\bS@>>>sU_P\cap S_s$ be the restriction of $\p:\bS@>>>S_s$. Now $U_P$ acts 
by conjugation on $suU_P\cap\bS$ and on $sU_P\cap S_s$ compatibly with $f$; moreover 
this action is transitive on $sU_P\cap S_s$ (see 19.3(a)). We have 
$f\i(s)=suU_Q\cap\bS$ hence we must only show that $H^i_c(f\i(s),K_0)=0$ for all $i$.
The Leray spectral sequence of $f$ is:
$$E_2^{p,q}=H^p_c(sU_P\cap S_s,\ch^qf_!K_0)\imp H^{p+q}_c(suU_P\cap\bS,K_0).$$
The last vector space is zero since $\ce$ is strongly cuspidal. Thus, $E_\iy^{p,q}=0$
for all $p,q$. Now $\ch^qf_!K_0$ is a $U_P$-equivariant local system on $sU_P\cap S_s$
and $sU_P\cap S_s\cong U_P/U_{P'}$ is an affine space. Hence $E_2^{p,q}=0$ for
$p\ne 2\dim U_P/U_{P'}$. This implies that $E_2^{p,q}=E_\iy^{p,q}$ for all $p,q$; it
follows that $E_2^{p,q}=0$ for all $p,q$ so that $\ch^qf_!K_0=0$ for all $q$. Taking 
the stalk at $s$ we see that $H^q_c(f\i(s),K_0)=0$ for all $q$ and (b) is proved.

We show:

(c) {\it Assume that $S={}^D\cz_{G^0}^0\boc$ where $\boc$ is a unipotent 
$G^0$-conjugacy class. Let $\cl\in\cs({}^D\cz_{G^0}^0)$ be a local system of rank $1$. 
If $\ce$ is strongly cuspidal then $\ce\ot(\cl\bxt\bbq)$ is strongly cuspidal.}
\nl
We have $IC(\bS,\ce\ot(\cl\bxt\bbq))=K_0\ot(\cl\bxt\bbq)$ since
$\bS\cong{}^D\cz_{G^0}^0\T\bboc$. Let $P$ be a parabolic of $G^0$ such that 
$P\ne G^0,S\sub N_GP$. Let $z\in{}^D\cz_{G^0}^0,u\in\bboc$. We know that 
$H^i_c(zuU_P\cap\bS,K_0)=0$. We must show that 
$H^i_c(zuU_P\cap\bS,K_0\ot(\cl\bxt\bbq))=0$. It is enough to show that 
$(\cl\bxt\bbq)|_{zuU_P\cap\bS}\cong\bbq$. This follows from the fact that 
$(\cl\bxt\bbq)|_{zuU_P\cap\bS}=\cl_z\ot\bbq$ where $\cl_z$ is the stalk of $\cl$ at 
$z$.

This argument shows also that, in the setup of (c):

(d) {\it  If $\ce$ is clean then $\ce\ot(\cl\bxt\bbq)$ is clean.}

\subhead 23.5\endsubhead
Let $S$ be an isolated stratum of $G$. Let $D$ be the connected component of $G$ that 
contains $S$. Assume that there exists a non-zero cuspidal local system 
$\ce_0\in\cs(S)$. Let $\ce',\ce''\in\cs(S)$ be local systems such that the local system
$\ce'\ot\ce''$ has no direct summand isomorphic to $\bbq$. We show that

(a) $H^i_c(S,\ce'\ot\ce'')=0$ {\it for all $i$.}
\nl
It is enough to show that, if $\ce\in\cs(S)$ is irreducible and $\ce\not\cong\bbq$ then
$H^i_c(S,\ce)=0$ for all $i$. Let $H={}^D\cz_{G^0}^0\T G^0$. We can find an integer 
$n\ge 1$ invertible in $\kk$ such that $\ce$ is equivariant for the transitive 
$H$-action $(z,x):g\m xz^ngx\i$ on $S$. Let $y\in S$ and let $H_y$ be the stabilizer of
$y$ in $H$ for this action. Let $\tS=H/H_y^0$. Define $f:\tS@>>>S$ by 
$f(z,x)=xz^nyx\i$ (a principal $H_y/H_y^0$-covering where $H_y$ acts on $\tS$ by right 
multiplication). Now $\ce$ is a direct summand of the local system $f_!\bbq$. It is 
enough to show that $H^i_c(S,f_!\bbq)=H^i_c(S,\bbq)$ for all $i$ or equivalently that 
$H_y/H_y^0$ acts trivially on $H^i_c(\tS,\bbq)$. Let 
$\bH={}^D\cz_{G^0}^0\T(G^0/{}^D\cz_{G^0}^0)$. Define $f':\bH@>>>\tS$ by 
$(z,x{}^D\cz_{G^0}^0)\m(z,x)H_y^0$. Now $f'$ is a fibration with fibres isomorphic to
$H_y^0/(\{1\}\T{}^D\cz_{G^0}^0)\cong Z_G(y)^0/{}^D\cz_{G^0}^0$ which by 10.2 is 
isomorphic to an affine space of dimension, say, $a$. Hence we have
$H^i_c(\tS,\bbq)=H^{i+2a}_c(\bH,\bbq)$. Also the $H_y/H_y^0$-action on $\tS$ is 
compatible under $f'$ with the $H_y/(\{1\}\T{}^D\cz_{G^0}^0)$-action on $\bH$ by right 
multiplication. It is then enough to show that $H_y/(\{1\}\T{}^D\cz_{G^0}^0)$ acts 
trivially on $H^{i+2a}_c(\bH,\bbq)$. This follows from the fact that the 
$H_y/(\{1\}\T{}^D\cz_{G^0}^0)$-action on $\bH$ is the restriction of an action of the 
connected group $H/(\{1\}\T{}^D\cz_{G^0}^0)$ and a connected group must act trivially 
in cohomology. This proves (a).

\subhead 23.6\endsubhead
Let $(L',S'),(L'',S'')\in\AA$. Assume that $\ce'\in\cs(S')$ and $\ce''\in\cs(S'')$ are 
strongly cuspidal relative to $N_GL',N_GL''$ respectively. Let $K'_0=IC(\bS',\ce')$,
$K''_0=IC(\bS'',\ce'')$. We regard $K'_0$ (resp. $K''_0$) as a complex on $N_GL'$ 
(resp. $N_GL''$) zero outside $\bS'$ (resp. $\bS''$). To $L',S',\ce'$ (resp.
$L'',S'',\ce''$) we attach $\fK'\in\cd(\bY_{L',S'})$ (resp. 
$\fK''\in\cd(\bY_{L'',S''})$) in the same way as $\fK=IC(\bY_{L,S},\p_!\tce)$ was
attached to $L,S,\ce$ in 5.6. We regard $\fK',\fK''$ as complexes on $G$, zero outside 
$\bY_{L',S'},\bY_{L'',S''}$ respectively. 

\proclaim{Proposition 23.7}Assume that for any $n\in G^0$ such that $n\i L'n=L''$, 
hence $n\i N_GL'n=N_GL''$, we have $H^i_c(N_GL',K'_0\ot\Ad(n\i)^*K''_0)=0$ for all $i$.
(This condition is automatically verified if $L',L''$ are not $G^0$-conjugate.) Then\lb
$H^i_c(G,\fK'\ot\fK'')=0$ for all $i$.
\endproclaim
The proof is quite similar to (but simpler than) that in 7.8. Let $P'$ (resp. $P''$) be
a parabolic of $G^0$ with Levi $L'$ (resp. $L''$) such that $S'\sub N_GP'$ (resp. 
$S''\sub N_GP''$). Let $X',X'',K'\in\cd(X'),K''\in\cd(X''),\fZ$ be as in 7.4. We may 
regard $\fZ$ as a subvariety of $X'\T X''$ via the imbedding 
$(g,x'P',x''P'')\m((g,x'P'),(g,x''P''))$. The inverse image of 
$K'\bxt K''\in\cd(X'\T X'')$ under this imbedding is a complex $\tK\in\cd(\fZ)$. Using 
a description of $\fK',\fK''$ like in 5.7 we see that $\fK'\ot\fK''=(pr_1)_!\tK$ where 
$pr_1:\fZ@>>>G$ is the first projection. It follows that
$H^i_c(G,\fK'\ot\fK'')=H^i_c(\fZ,\tK)$ for all $i$. Hence it is enough to show that 
$H^i_c(\fZ,\tK)=0$ for all $i$. For any $G^0$-orbit $E$ on $G^0/P'\T G^0/P''$ let 
${}^E\fZ=\{(g,x'P',x''P'')\in\fZ;(x'P',x''P'')\in E\}$. Using the partition of $\fZ$ 
into the finitely many locally closed subvarieties ${}^E\fZ$ we see that it is enough 
to show that $H^i_c({}^E\fZ,\tK)=0$ for all $i$ and any $E$ as above. Using the
spectral sequence of the fibration $pr_{23}:{}^E\fZ@>>>E$, 
$(g,x'P',x''P'')\m(x'P',x''P'')$ we see that it is enough to show that for any 
$(x',x'')\in G^0\T G^0$ such that $(x'P',x''P'')\in E$ we have 

$H^i_c(\cv,\tK)=0$ for all $i$
\nl
where $\cv$ is the fibre of $pr_{23}$ at $(x'P',x''P'')$. We may identify

$\cv=\{g\in G;x'{}\i gx'\in\bS'U_{P'},x''{}\i gx''\in\bS''U_{P''}\}$.
\nl
Define $j:\cv@>>>\bS'\T\bS''$ by 
$$j(g)=(\bS'-\text{component of }x'{}\i gx',\bS''-\text{component of }x''{}\i gx'').$$
Then $\tK|_\cv$ may be identified with $j^*(K'_0\bxt K''_0)$ and we must show that

$H^i_c(\cv,j^*(K'_0\bxt K''_0))=0$ for all $i$.
\nl
Let $Q',Q'',M',M'',\Si',\Si'',\cf',\cf''$ be as in 7.8. Let $\hK'_0=IC(\bSi',\cf')$,
$\hK''_0=IC(\bSi'',\cf'')$. As in 7.8, $\cv$ is fibred over
$$\cv_1=\{(u'',u',z)\in(M'\cap U_{Q''})\T(M''\cap U_{Q'})\T(\tM'\cap\tM'');
zu''\in\bSi',zu'\in\bSi''\}$$
with all fibres isomorphic to $U_{Q'}\cap U_{Q''}$. Since $U_{Q'}\cap U_{Q''}$ is an
affine space we see that it is enough to show that
$H^i_c(\cv_1,\bj^*(\hK'_0\bxt\hK''_0))=0$ for all $i$ where 
$\bj:\cv_1@>>>\bSi'\T\bSi''$ is defined by $\bj(u'',u',z)=(zu'',zu')$.

Assume first that $Q',Q''$ have no common Levi. Let $\tp_3:\cv_1@>>>\tM'\cap\tM''$ be 
the third projection. It is enough to show that for any $z\in\tM'\cap\tM''$ we have 
$H^i_c(\tp_3\i(z),\bj^*(\hK'_0\bxt\hK''_0))=0$ for all $i$. Now $\tp_3\i(z)$ is a 
product $R'\T R''$ where $R'$ (resp. $R''$) is the set of all elements in 
$(\tM'\cap N_GQ'')\cap\bSi'$ (resp. $(\tM''\cap N_GQ')\cap\bSi''$) whose image under 
$\tM'\cap N_GQ''@>>>\tM'\cap\tM''$ (resp. $\tM''\cap N_GQ'@>>>\tM'\cap\tM''$) is equal 
to $z$. We are reduced to showing that 
$H^{i'}_c(R',\hK'_0)\ot H^{i''}_c(R'',\hK''_0)=0$ for all $i',i''$. Since $Q',Q''$ have
no common Levi, we see that either $M'\cap Q''$ is a proper parabolic of $M'$ or 
$M''\cap Q'$ is a proper parabolic of $M''$. In the first case we have 
$H^{i'}_c(R',\hK'_0)=0$ for all $i'$ since $\cf'$ is a strongly cuspidal local system 
in $\cs(\Si')$. In the second case we have $H^{i''}_c(R'',\hK''_0)=0$ for all $i''$ 
since $\cf''$ is a strongly cuspidal local system in $\cs(\Si'')$. Thus the desired 
vanishing result holds in our case.

Assume next that $Q',Q''$ have a common Levi. Then $M'=M''$, $M'\cap U_{Q''}=\{1\}$,
$M''\cap U_{Q'}=\{1\}$ and we may identify $\cv_1$ with $\bSi'\cap\bSi''$. Hence we 
must only show that $H^i_c(\bSi'\cap\bSi'',\hK'_0\ot\hK''_0))=0$ for all $i$. As in 7.8
we can find $v'\in U_{P'},v''\in U_{P''}$ such that, setting
$n=v'{}\i x'{}\i x''v''\in G^0$, we have $n\i L'n=L''$. Then the desired vanishing 
result is equivalent to $H^i_c(N_GL',K'_0\ot\Ad(n\i)^*K''_0)=0$ for all $i$, which is 
part of our assumptions. The proposition is proved.

\head 24. Orthogonality\endhead
\subhead 24.1\endsubhead
In this section we assume that $\kk$ is an algebraic closure of a finite field $\FF_q$
and that $G$ has a fixed $\FF_q$-rational structure with Frobenius map $F:G@>>>G$. 

\subhead 24.2\endsubhead
Let $L$ be a Levi of a parabolic of $G^0$. Let $\d$ be a connected component of $N_GL$
such that $\d\sub N_G^\bul L$. Assume that $F(L)=L,F(\d)=\d$ and that the $F$-stable 
torus $\cz={}^\d\cz_L^0$ is split over $\FF_q$ that is, $F(z)=z^q$ for all $z\in\cz$. 
We show:

(a) {\it there exists a parabolic $P$ of $G^0$ with Levi $L$ such that $F(P)=P$.}
\nl
We can find $\c\in\Hom(\kk^*,\cz)$ such that $Z_{G^0}(\c(\kk^*))=Z_{G^0}(\cz)$. Let 
$g\in\d$. Define $\c',\c''\in\Hom(\kk^*,\cz)$ by $\c'(a)=g\c(a)g\i$,
$\c''(a)=F(\c(a^{q\i}))$. To $\c,\c',\c''$ we attach parabolics $P_\c,P_{\c'},P_{\c''}$
of $G^0$ as in 1.16. Let $P=P_\c$. From the definitions we have $P_{\c'}=gPg\i$,
$P_{\c''}=F(P)$. Since $\c(\kk^*)\sub Z_{G^0}(g)$ we have $\c'=\c$. Since $\cz$ is 
$\FF_q$-split, we have $\c''(a)=\c(a^{q\i})^q=\c(a)$ for any $a$ hence $\c''=\c$. Thus,
$gPg\i=P,F(P)=P$. It follows that $\d\sub N_GP$. Now $Z_{G^0}(\c(\kk^*))$ is a Levi of 
$P$ and $Z_{G^0}(\cz)=L$ by 1.10(a). Hence $L$ is a Levi of $P$. This proves (a).

\subhead 24.3\endsubhead
Let $(L',S'),(L'',S'')\in\AA$. Assume that $\ce'\in\cs(S')$ and $\ce''\in\cs(S'')$ are 
strongly cuspidal relative to $N_GL',N_GL''$ respectively. Assume that
$F(L')=L',F(S')=S',F(L'')=L'',F(S'')=S''$ and that we are given isomorphisms 
$\e':F^*\ce'@>\si>>\ce',\e'':F^*\ce''@>\si>>\ce''$. Let $\fK',\fK''$ be as in 23.6 and 
let $\ph':F^*\fK'@>\si>>\fK'$, $\ph'':F^*\fK''@>\si>>\fK''$ be the isomorphisms 
induced by $\e',\e''$. Let $\d'$ (rep. $\d''$) be the connected component of $N_GL'$ 
(resp. $N_GL''$) that contains $S'$ (resp. $S''$). Assume that either $L',L''$ are not 
$G^0$-conjugate or that $\ce'$ and $\ce''$ are clean (relative to $N_GL',N_GL''$). Let 
$\Th=\{n\in G^{0F};n\i L'n=L'',n\i S'n=S''\}$. With these assumptions we state the 
following result.

\proclaim{Lemma 24.4}$$\align&|G^{0F}|\i\sum_{g\in G^F}\c_{\fK',\ph'}(g)
\c_{\fK'',\ph''}(g)\\&=|L'{}^F|\i|L''{}^F|\i\sum_{n\in\Th}\sum_{y\in S'{}^F}
\c_{\ce',\e'}(y)\c_{\ce'',\e''}(n\i yn).\tag a\endalign$$
\endproclaim
The proof is given in 24.7-24.12.

\subhead 24.5\endsubhead
Let $(L',S'),(L'',S'')\in\AA$. Let $\d'$ (resp. $\d''$) be the connected component of
$N_GL'$ (resp. $N_GL''$) that contains $S'$ (resp. $S''$). Assume that 
$S'={}^{\d'}\cz_{L'}^0\T\boc',S''={}^{\d''}\cz_{L''}^0\T\boc''$ where $\boc'$ is a 
unipotent $L'$-conjugacy class and $\boc''$ is a unipotent $L''$-conjugacy class. Let 
$\cf'$ (resp. $\cf''$) be an $L'$- (resp. $L''$-) equivariant local system on $\boc'$ 
(resp. $\boc''$) such that $\bbq\bxt\cf'$ (resp. $\bbq\bxt\cf''$) is a strongly 
cuspidal local system relative to $N_GL'$ (resp. $N_GL''$). Assume that $F(L')=L'$,
$F(\boc')=\boc'$, $F(L'')=L''$, $F(\boc'')=\boc''$ and that we are given isomorphisms 
$\bae':F^*\cf'@>\si>>\cf'$, $\bae'':F^*\cf''@>\si>>\cf''$. Assume that either $L',L''$ 
are not $G^0$-conjugate or that $\bbq\bxt\cf'$ and $\bbq\bxt\cf''$ are clean (relative 
to $N_GL',N_GL''$). Let $\Th$ be as in 24.3. With these assumptions we state the 
following result.

\proclaim{Lemma 24.6}
$$\align&|G^{0F}|\i\sum_{u\in G^F_{un}}Q_{L',G,\boc',\cf',\bae'}(u)
Q_{L'',G,\boc'',\cf'',\bae''}(u)\\&=|L'{}^F|\i|L''{}^F|\i\sum_{n\in\Th}
\sum_{y\in\boc'{}^F}\c_{\cf',\bae'}(y)\c_{\cf'',\bae''}(n\i yn).\tag a\endalign$$
\endproclaim
The proof is given (together with that of Lemma 24.4) in 24.7-24.12.

\subhead 24.7\endsubhead
We prove Lemma 24.4 assuming that $L',L''$ are not $G^0$-conjugate. Let
$\Ph:H^i_c(G,\fK'\ot\fK'')@>\si>>H^i_c(G,\fK'\ot\fK'')$ be the composition

$H^i_c(G,\fK'\ot\fK'')@>\si>>H^i_c(G,F^*\fK'\ot F^*\fK'')@>\si>>H^i_c(G,\fK'\ot\fK'')$ 
\nl
(the first map is induced by $F:G@>>>G$, the second map is induced by $\ph'\ot\ph''$). 
By the trace formula for Frobenius maps, the left hand side of 24.4(a) is equal to 
$\sum_i(-1)^i\tr(\Ph,H^i_c(G,\fK'\ot\fK''))$. This is zero since
$H^i_c(G,\fK'\ot\fK'')=0$ for all $i$, by 23.7. The right hand side of 24.4(a) is also
zero, by our assumption. Thus Lemma 24.4 is proved in the present case.

\subhead 24.8\endsubhead
We prove Lemma 24.4 under the following assumptions:

$S'={}^{\d'}\cz_{L'}^0\boc',S''={}^{\d''}\cz_{L''}^0\boc''$ where 
$\d',\d',\boc',\boc''$ are as in 24.5;

for some/any $c'\in\boc'$, $\ce'|_{{}^{\d'}\cz_{L'}^0c'}$ has no direct summand 
isomorphic to $\bbq$;

for some/any $c''\in\boc''$, $\ce''|_{{}^{\d''}\cz_{L''}^0c''}$ is isomorphic to 
$\bbq^N$ for some $N$.
\nl
By 24.7 we may assume that $L,L'$ are $G^0$-conjugate. Then $\ce',\ce''$ are clean. We
show that the left hand side of 24.4(a) is zero. As in 24.7 it is enough to show that 
$H^i_c(G,\fK'\ot\fK'')=0$ for all $i$. Using 23.7, it is enough to show that for any 
$n\in G^0$ such that $n\i L'n=L''$ we have $H^i_c(N_GL',K'_0\ot\Ad(n\i)^*K''_0)=0$ for 
all $i$. The last equality is clear if $n\i S'n\ne S''$ (in this case 
$n\i S'n\cap S''=\em$ and we have $K'_0\ot\Ad(n\i)^*K''_0=0$ since $\ce',\ce''$ are 
clean). Assume now that $n\i S'n=S''$. Then $K'_0\ot\Ad(n\i)^*K''_0$ is the local 
system $\ce'\ot\Ad(n\i)^*\ce''$ on $S'$ extended by $0$ on $N_GL'-S'$. Hence it is 
enough to show that $H^i_c(S',\ce'\ot\Ad(n\i)^*\ce'')=0$ for all $i$. This follows from
23.5(a) since, by our assumption, the local system $\ce'\ot\Ad(n\i)^*\ce''$ has no 
direct summand isomorphic to $\bbq$. Next we show that the right hand side of 24.4(a) 
is zero. It is enough to show that for any $n\in\Th$ the sum 
$\sum_{y\in S'{}^F}\c_{\ce',\e'}(y)\c_{\ce'',\e''}(n\i yn)$ is zero. By the trace 
formula for Frobenius maps, this sum is an alternating sum of traces of the Frobenius 
map on \lb $H^i_c(S',\ce'\ot\Ad(n\i)^*\ce'')$. As we have seen above, this last vector 
space is zero. Thus Lemma 24.4 is proved in the present case.

\subhead 24.9\endsubhead
We prove Lemma 24.4 under the following assumption: there exist parabolics $P',P''$ of 
$G^0$ with Levi $L',L''$ respectively such that $S'\sub N_GP'S''\sub N_GP''$,
$F(P')=P',F(P'')=P''$.

By 24.7 we may assume that $L,L'$ are $G^0$-conjugate. Then $\ce',\ce''$ are clean.
Define $X'_{S'},\bce'$ (resp. $X''_{S''},\bce''$) in terms of $(L',P',S',\ce')$ (resp. 
$(L'',P'',S'',\ce'')$) in the same way as $X_S,\bce$ were defined in terms of 
$(L,P,S,\ce)$ in 5.6. Using 5.7 and the cleanness of $\ce',\ce''$ we see that
$\fK'=f'_!\bce',\fK''=f''_!\bce''$ where $f':X'_{S'}@>>>G,f'':X''_{S''}@>>>G$ are given
by the first projection. Note that $X'_{S'},X''_{S''},f',f''$ are naturally defined 
over $\FF_q$ and there are obvious isomorphisms $F^*\bce'@>\si>>\bce'$,
$F^*\bce''@>\si>>\bce''$ induced by $\e',\e''$. It follows that
$$\c_{\fK',\ph'}(g)=\sum_{x'P'{}^F\in G^{0F}/P'{}^F;x'{}\i gx'\in S'U_{P'}}
\c_{\ce',\e'}(\p'(x'{}\i gx')),$$
$$\c_{\fK'',\ph''}(g)=\sum_{x''P''{}^F\in G^{0F}/P''{}^F;x''{}\i gx''\in S''U_{P''}}
\c_{\ce'',\e''}(\p''(x''{}\i gx'')),$$
where $\p':(S'U_{P'})^F@>>>S'{}^F,\p'':(S''U_{P''})^F@>>>S''{}^F$ are the obvious 
projections (see 1.26). Hence the left hand side of 24.4(a) is
$$|G^{0F}|\i|P'{}^F|\i|P''{}^F|\i\sum_{x',x''\in G^{0F}}h(x',x'')\tag a$$
where
$$h(x',x'')=\sum_{g\in V^F}\c_{\ce',\e'}(\p'(x'{}\i gx'))\c_{\ce'',\e''}
(\p''(x''{}\i gx''))$$
and $V=\{g\in G;x'{}\i gx'\in S'U_{P'},x''{}\i gx''\in S''U_{P''}\}$.

Assume first that $x'P'x'{}\i,x''P''x''{}\i$ have no common Levi. In this case we show
that $h(x',x'')=0$. By the trace formula for Frobenius maps, $h(x',x'')$ is equal to an
alternating sum of traces of Frobenius on $H^i_c(\cv,\tK)$ (notation as in the proof of
23.7). But by the proof of 23.7, in our case we have $H^i_c(\cv,\tK)=0$ for all $i$.
(The relevant part of the proof of 23.7 does not make use of the assumptions in the 
first sentence of 23.7; it only uses the strong cuspidality of $\ce',\ce''$.) Thus 
$h(x',x'')=0$, as desired. 

Assume now that $x',x''\in G^{0F}$ are such that $Q'=x'P'x'{}\i,Q''=x''P''x''{}\i$ have
a common Levi $M'=M''$ (we may assume that $M'=M''$ is $F$-stable). Then we can find 
$v'\in U_{P'}^F,v''\in U_{P''}^F$ such that, setting 
$n=v'{}\i x'{}\i x''v''\in G^{0F}$, we have $n\i L'n=L''$. Also, $n$ is uniquely 
determined by $x',x''$. Let $\Si',\Si''$ be strata of $N_GM'=N_GM''$ as in 7.8. As in 
7.8, we have a natural map $V@>>>V_1=\Si'\cap\Si''$ with fibres isomorphic to 
$U_{Q'}\cap U_{Q''}\cong U_{P'}\cap U_{nP''n\i}$. In particular, $V=\em$ unless 
$\Si'\cap\Si''\ne\em$ or equivalently, $\Si'=\Si''$ or equivalently, $n\i S'n=S''$. If 
this last condition is satisfied, we have 
$$h(x',x'')=|U_{P'}^F\cap U_{nP''n\i}^F| 
\sum_{y\in S'{}^F}\c_{\ce',\e'}(y)\c_{\ce'',\e''}(n\i yn).$$
It is then enough to show that
$$\align&\sh((x',x'')\in G^{0F}\T G^{0F};x'{}\i x''\in U_{P'}^FnU_{P''}^F)\\&=
|G^{0F}||P'{}^F||P''{}^F||L'{}^F|\i|L''{}^F|\i |U_{P'}^F\cap U_{nP''n\i}^F|.\endalign$$
This is immediate. Thus Lemma 24.4 is proved in the present case.

\subhead 24.10\endsubhead
We show that Lemma 24.4 holds for $G$ under the assumption that Lemma 24.6 holds when 
$G$ is replaced by $Z_G(s)$ where $s$ is any semisimple element of $G^F$. We evaluate 
the left hand side of 24.4(a) using the "character formula" 16.14. We have
$$\align&|G^{0F}|\i\sum_{g\in G^F}\c_{\fK',\ph'}(g)\c_{\fK'',\ph''}(g)=
|G^{0F}|\i\sum_{s\in G^F \text{semis.}}\sum \Sb x',x''\in G^{0F}\\
x'{}\i sx'\in S'_s\\x''{}\i sx''\in S''_s\endSb       \\&
|Z_G(s)^{0F}|^{-2}|L'_{x'}{}^F||L'{}^F|\i|L''_{x''}{}^F||L''{}^F|\i 
f(s,x',x'',\dd',\dd'')\endalign$$
where

$L'_{x'}=x'L'x'{}\i\cap Z_G(s)^0$, $L''_{x''}=x''L''x''{}\i\cap Z_G(s)^0$,

$\dd'$ runs over the set of $F$-stable $L'_{x'}$-conjugacy classes contained in

$\{v\in Z_G(s);v\text{ unipotent, }x'{}\i svx'\in S'\}$, 

$\dd''$ runs over the set of $F$-stable $L''_{x''}$-conjugacy classes contained in

$\{v\in Z_G(s);v\text{ unipotent, }x''{}\i svx''\in S''\}$, and
$$f(s,x',x'',\dd',\dd'')=\sum\Sb u\in Z_G(s)^F\\ u\text{unip.}\endSb
Q_{L'_{x'},Z_G(s),\dd',\cf'_{x'},\e'_{x'}}(u)
Q_{L''_{x''},Z_G(s),\dd'',\cf''_{x''},\e''_{x''}}(u).$$
Here $\cf'_{x'}$ is the inverse image of $\ce'$ under $\dd'@>>>S',v\m x'{}\i svx'$, 
$\cf''_{x''}$ is the inverse image of $\ce''$ under $\dd''@>>>S'',v\m x''{}\i svx''$,
and $\e'_{x'}:F^*\cf'_{x'}@>\si>>\cf'_{x'},\e''_{x''}:F^*\cf''_{x''}@>\si>>\cf''_{x''}$
are induced by $\e',\e''$.

Using our assumption we have
$$\align&f(s,x',x'',\dd',\dd'')=|Z_G(s)^{0F}||L'_{x'}{}^F|\i|L''_{x''}{}^F|\i
\sum\Sb n\in Z_G(s)^{0F}\\ n\i L'_{x'}n=L''_{x''}\\n\i\dd' n=\dd''\endSb\\& \T
\sum_{v\in\dd'{}^F}\c_{\cf'_{x'},\e'_{x'}}(v)\c_{\cf''_{x''},\e''_{x''}}(n\i vn).\tag a
\endalign$$
To be able to apply our assumption, we use 23.4(a),(b). We also use the following fact.

{\it For $s,x',x''$ as above and for $n\in Z_G(s)^0$, the condition 
$n\i L'_{x'}n=L''_{x''}$ is equivalent to $n\i x'Lx'{}\i n=x''Lx''$.}
\nl
Indeed, since $S'$ is an isolated stratum of $N_GL'$ and $x'{}\i sx'\in S_s$, we see 
using 18.2 that $x'{}\i sx'$ is isolated in $N_GL'$. Hence $s$ is isolated in 
$N_G(x'L'x'{}\i)$. It follows that $s$ is isolated in $N_G(n\i x'L'x'{}\i n)$. 
Similarly, $s$ is isolated in $N_G(x''L''x''{}\i)$. By the injectivity of the map $a$ 
in 21.3 (for $s$ instead of $g$) we see that  
$$n\i x'L'x'{}\i n\cap Z(s)^0=x''L''x''{}\i\cap Z_G(s)^0\lra 
n\i x'L'x'{}\i n=x''L''x''{}\i,$$
that is $n\i L'_{x'}n=L''_{x''}\lra  n\i x'L'x'{}\i n=x''L''x''{}\i$, as required. We 
have
$$\c_{\cf'_{x'},\e'_{x'}}(v)=\c_{\ce',\e'}(x'{}\i svx'),\quad
\c_{\cf''_{x''},\e''_{x''}}(n\i vn)=\c_{\ce'',\e''}(x''{}\i sn\i vnx'').$$
Note also that the condition $n\i\dd' n=\dd''$ implies 
$n\i x'S'x'{}\i n=x''S''x''{}\i$. We see that
$$\align&|G^{0F}|\i\sum_{g\in G^F}\c_{\fK',\ph'}(g)\c_{\fK'',\ph''}(g)\\&=
|G^{0F}|\i|L'{}^F|\i|L''{}^F|\i\sum\Sb s\in G^F \text{semis.}\\
x',x''\in G^{0F},\dd'\\x'{}\i sx'\in S'_s\\x''{}\i sx''\in S''_s\endSb|Z_G(s)^{0F}|\i
\\& \sum\Sb n\in Z_G(s)^{0F}\\n\i x'L'x'{}\i n=x''L''x''{}\i\\
n\i x'S'x'{}\i n=x''S''x''{}\i\endSb\sum_{v\in\dd'{}^F}
\c_{\ce',\e'}(x'{}\i svx')\c_{\ce'',\e''}(x''{}\i sn\i vnx'').\endalign$$
We now make the change of variable $(x',x'',n)\m(x',n,n'),n'{}\i=x''{}\i n\i x'$. The 
condition $n\i x'L'x'{}\i n=x''L''x''{}\i$ becomes $n'{}\i L'n'=L''$; the condition
$n\i x'S'x'{}\i n=x''S''x''{}\i$ becomes $n'{}\i S'n'=S''$ (thus, $n'\in\Th$). The
condition $x''{}\i sx''\in S''_s$ becomes $n'{}\i x'{}\i nsn\i x'n'\in n'{}\i S'_sn'$ 
that is $x'{}\i sx'\in S'_s$. Our sum becomes
$$\align&|G^{0F}|\i|L'{}^F|\i|L''{}^F|\i 
\sum\Sb s\in G^F \text{semis.}\\x'\in G^{0F}\\x'{}\i sx'\in S'_s\endSb|Z_G(s)^{0F}|\i
\\&\sum\Sb n\in Z_G(s)^{0F}\\n'\in\Th\endSb\sum\Sb v\in Z_G(s)^F; v\text{ unip.}\\
x'{}\i svx'\in S'\endSb\c_{\ce',\e'}(x'{}\i svx')
\c_{\ce'',\e''}(n'{}\i x'{}\i sv x'n').\endalign$$
We make the change of variable $(s,x',v)\m(s',x',v')$ where $s'=x'{}\i sx'\in S'_s$, 
$v'=x'{}\i vx'\in Z_G(s')\cap s'{}\i S'$. Our sum becomes
$$\align&|L'{}^F|\i|L''{}^F|\i\sum\Sb s'\in S'_s{}^F\\n'\in\Th\endSb
\sum\Sb v'\in Z_G(s')^F;v'\text{ unip.}\\v'\in s'{}\i S'\endSb\c_{\ce',\e'}(s'v')
\c_{\ce'',\e''}(n'{}\i s'v'n')\\&=|L'{}^F|\i|L''{}^F|\i\sum_{n'\in\Th}
\sum_{y\in S'{}^F}\c_{\ce',\e'}(y)\c_{\ce'',\e''}(n'{}\i yn'),\endalign$$
as required.

\subhead 24.11\endsubhead
We show that Lemma 24.6 holds for $G$ under the assumption that Lemma 24.6 holds when 
$G$ is replaced by $Z_G(s)$ where $s$ is any semisimple element of $G^F$ such that 
$\dim Z_G(s)<\dim G$ (that is, $s\n Z_G(G^0)$). 

Let $\cl'\in\cs({}^{\d'}\cz_{L'}^0)$, $\cl''\in\cs({}^{\d''}\cz_{L''}^0)$ be local
systems of rank $1$ with given isomorphisms $\io':F^*\cl'@>\si>>\cl'$,
$\io'':F^*\cl''@>\si>>\cl''$. Let $\ce'=\cl'\bxt\cf'\in\cs(S')$,
$\ce''=\cl''\bxt\cf''\in\cs(S'')$. Let $\e'=\io'\bxt\bae':F^*\ce'@>\si>>\ce'$, 
$\e''=\io''\bxt\bae'':F^*\ce''@>\si>>\ce''$. Note that $\ce',\ce''$ are strongly
cuspidal by 23.4(c) and that, if $L',L''$ are not $G^0$-conjugate then $\ce',\ce''$ are
clean by 23.4(d). For this $\ce',\ce''$ we can still try to carry out the 
argument in 24.10 but now we can only use 24.6(a) for $s$ such that  $s\n Z_G(G^0)$. We
obtain:
$$\align&|G^{0F}|\i\sum_{g\in G^F}\c_{\fK',\ph'}(g)\c_{\fK'',\ph''}(g)\\&
-|L'{}^F|\i|L''{}^F|\i\sum_{n\in\Th}\sum_{y\in S'{}^F}
\c_{\ce',\e'}(y)\c_{\ce'',\e''}(n\i yn)\\&
=|G^{0F}|\i\sum\Sb s\in Z_G(G^0)^F \text{semis.}\endSb
\sum\Sb x',x''\in G^{0F},\dd',\dd''\\x'{}\i sx'\in S'_s\\x''{}\i sx''\in S''_s\endSb\\&
|Z_G(s)^{0F}|^{-2}|L'_{x'}{}^F||L'{}^F|\i|L''_{x''}{}^F||L''{}^F|\i \\&
(\sum_{u\in Z_G(s)^F;\text{unip.}}Q_{L'_{x'},Z_G(s),\dd',\cf'_{x'},\e'_{x'}}(u)
Q_{L''_{x''},Z_G(s),\dd'',\cf''_{x''},\e''_{x''}}(u)\\&-|Z_G(s)^{0F}||L'_{x'}{}^F|\i
|L''_{x''}{}^F|\i\sum\Sb n\in Z_G(s)^{0F}\\n\i L'_{x'}n=L''_{x''}\\n\i\dd'n=\dd''\endSb
\\&\sum_{v\in\dd'{}^F}\c_{\cf'_{x'},\e'_{x'}}(v)\c_{\cf''_{x''},\e''_{x''}}(n\i vn))
\tag a\endalign$$
(notation of 24.10.) For each $s,x',x''$ in the right hand side of (a) we have 
$s=x'{}\i sx'\in{}^{\d'}\cz_{L'}^0$. Similarly, $s\in{}^{\d''}\cz_{L''}^0$. In 
particular we have $s\in G^0$ hence $s\in\cz_{G^0}$ and $Z_G(s)^0=G^0$, 
$L'_{x'}=x'L'x'{}\i$, $L''_{x''}=x''L''x''{}\i$. Also we necessarily have 
$\dd'=x'\boc'x'{}\i$, $\dd''=x''\boc''x''{}\i$. We see that the right hand side of (a) 
is
$$\align&|G^{0F}|\i\sum_{s\in\cz_{G^0}\cap{}^{\d'}\cz_{L'}^0\cap{}^{\d''}\cz_{L''}^0}
\sum_{x',x''\in G^{0F}}\\&|G^{0F}|^{-2}(\sum_{u\in Z_G(s)^F;\text{unip.}}
Q_{L',Z_G(s),\boc',\cf',\bae'}(u)Q_{L'',Z_G(s),\boc'',\cf'',\bae''}(u)\\&
-|G^{0F}||L'{}^F|\i|L''{}^F|\i
\sum\Sb n\in G^{0F}\\n\i x'L'x'{}\i n=x''L''x''{}\i\\n\i x'\boc' x'{}\i n=\boc''\endSb
\sum_{v\in\boc'{}^F}\c_{\cf',\bae'}(v)\c_{\cf'',\bae''}(x''{}\i n\i x' vn))\endalign$$
or equivalently
$$\align&|G^{0F}|\i\sum_{s\in\cz_{G^0}\cap{}^{\d'}\cz_{L'}^0\cap{}^{\d''}\cz_{L''}^0}
\\&(\sum_{u\in Z_G(s)^F;\text{unip.}}Q_{L',Z_G(s),\boc',\cf',\bae'}(u)
Q_{L'',Z_G(s),\boc'',\cf'',\bae''}(u)\\&-|G^{0F}||L'{}^F|\i|L''{}^F|\i\sum_{n'\in\Th}
\sum_{v\in\boc'{}^F}\c_{\cf',\bae'}(v)\c_{\cf'',\bae''}(n'{}\i vn')).\tag b\endalign$$
Let $D'$ (resp. $D''$) be the connected component of $G$ that contains $\d'$ (resp. 
$\d''$) For each $s$ in the sum we have $D'\sub Z_G(s),D''\sub Z_G(s)$; moreover
$Q_{L',Z_G(s),\boc',\cf',\bae'}(u)=0$ unless $u\in D'$ and 
$Q_{L'',Z_G(s),\boc'',\cf'',\bae''}(u)=0$ unless $u\in D''$. We see that
$$\align&\sum_{u\in Z_G(s)^F;\text{unip.}}Q_{L',Z_G(s),\boc',\cf',\bae'}(u)
Q_{L'',Z_G(s),\boc'',\cf'',\bae''}(u)\\&=\sum_{u\in G^F;\text{unip.}} 
Q_{L',G,\boc',\cf',\bae'}(u)Q_{L'',G,\boc'',\cf'',\bae''}(u).\endalign$$
Hence (a) becomes
$$\align&|G^{0F}|\i\sum_{g\in G^F}\c_{\fK',\ph'}(g)\c_{\fK'',\ph''}(g)\\&
-|L'{}^F|\i|L''{}^F|\i\sum_{n\in\Th}\sum_{y\in S'{}^F}
\c_{\ce',\e'}(y)\c_{\ce'',\e''}(n\i yn)\\&
=|G^{0F}|\i\sh(\cz_{G^0}\cap{}^{\d'}\cz_{L'}^0\cap{}^{\d''}\cz_{L''}^0)
\\&(\sum_{u\in G^F;\text{unip.}}Q_{L',G,\boc',\cf',\bae'}(u)
Q_{L'',G,\boc'',\cf'',\bae''}(u)\\&-|G^{0F}||L'{}^F|\i|L''{}^F|\i\sum_{n'\in\Th}
\sum_{v\in\boc'{}^F}\c_{\cf',\bae'}(v)\c_{\cf'',\bae''}(n'{}\i vn')).\tag c\endalign$$
Hence to prove the equality 24.6(a) it is enough to show that the left hand side of (c)
is zero. In order to do so, we are free to choose $\cl',\cl'',\io',\io''$ in a
convenient way.

Assume first that ${}^{\d'}\cz_{L'}^{0F}\ne\{1\}$. Then we can find a non-trivial 
character $\th':{}^{\d'}\cz_{L'}^{0F}@>>>\bbq^*$ and $\cl',\io'$ as above such that 
$\c_{\cl',\io'}=\th'$. We have $\cl'\not\cong\bbq$. Let $\cl''=\bbq$ and take any 
$\io''$. With these choices, 24.8 shows that the left hand side of (c) is zero.

Assume next that ${}^{\d''}\cz_{L''}^{0F}\ne\{1\}$. Since $L',L''$ play a symmetrical
role we see as in the previous paragraph that the left hand side of (c) is zero.

Finally, assume that ${}^{\d'}\cz_{L'}^{0F}=\{1\}$ and 
${}^{\d''}\cz_{L''}^{0F}=\{1\}$. Then the $F$-stable tori ${}^{\d'}\cz_{L'}^0$,
${}^{\d''}\cz_{L''}^0$ are necessarily split over $\FF_q$ and $q=2$. Using 24.2(a) we 
can find parabolics $P',P''$ of $G^0$ with Levi $L',L''$ such that 
$F(P')=P',F(P'')=P''$. Then 24.9 shows that the left hand side of (c) is zero. Thus 
Lemma 24.6 is proved in the present case.

\subhead 24.12\endsubhead
Clearly, the arguments in 24.10, 24.11 provide an inductive proof of Lemmas 24.4, 24.6.

\subhead 24.13\endsubhead
In the setup of 24.3, assume that $\ce',\ce''$ are irreducible. Let 
$$\Th(\ce',\ce'')=\{n\in G^{0F};n\i L'n=L'',n\i S'n=S'',\Ad(n\i)^*\ce''\cong\che\ce'\}.
$$
If $n\in\Th(\ce',\ce'')$ the local system $\hce=\ce'\ot\Ad(n\i)^*\ce''$ is canonically 
of the form $\bbq\op\hce_1$ where $\hce_1\in\cs(S)$ has no direct summand $\bbq$. The 
isomorphisms $\e',\e''$ induce an isomorphism $F^*\hce@>\si>>\hce$ which respects the 
summand $\bbq$ and induces on it $\z(n)$ times the obvious isomorphism 
$F^*\bbq@>\si>>\bbq$. Here $\z(n)\in\bbq^*$ is well defined.

We show that, if $n\in\Th(\ce',\ce'')$ and $n_0\in G^{0F}$, $n_0\i L'n_0=L'$, 
$n_0\i S'n_0=S'$, $\Ad(n_0\i)^*\ce'\cong\ce'$, then $n_0n\in\Th(\ce',\ce'')$ satisfies

(a) $\z(n_0n)=\et(n_0\i)\z(n)$
\nl
where $\et(n_0\i)$ is as in 21.6.
Let $\a:\Ad(n_0\i)^*\ce'@>\si>>\ce'$ be an isomorphism. We have an isomorphism
$$\align&\Ad(n_0\i)^*(\ce'\ot\Ad(n\i)^*\ce'')\\&=\Ad(n_0\i)^*\ce'\ot\Ad(n\i n_0\i)^*
\ce''@>\a\ot 1>>\ce'\ot\Ad(n\i n_0\i)^*\ce''\endalign$$
which must carry the summand $\Ad(n_0\i)^*(\bbq)$ to the summand $\bbq$. Let 
$x\in S'{}^F$. Let $e'_i\in\ce'_{n_0\i xn_0},e''_i\in\ce''_{n\i n_0\i xn_0n}$ be such 
that $\sum_ie'_i\ot e''_i$ belongs to the fibre of the summand $\bbq$ of 
$\Ad(n_0\i)^*\ce'\ot\Ad(n\i n_0\i)^*\ce''$ at $x$. Then $\sum_i\a(e'_i)\ot e''_i$ 
belongs to the fibre of the summand $\bbq$ of $\ce'\ot\Ad(n\i n_0\i)^*\ce''$ at $x$. By
definition,
$$\sum_i\e'(e'_i)\ot\e''(e''_i)=\z(n)\sum_ie'_i\ot e''_i,
\sum_i\e'(\a(e'_i))\ot\e''(e''_i)=\z(n_0n)\sum_i\a(e'_i)\ot e''_i,$$
$$\e'(\a(e'_i))=\et(n_0\i)\a(\e'(e'_i)).$$
\nl
We deduce
$$\align&\sum_i\e'(\a(e'_i))\ot\e''(e''_i)\\&=\et(n_0\i)(\a\ot 1)(\sum_i\e'(e'_i))
\ot\e''(e''_i))\et(n_0\i)\z(n)(\a\ot 1)(\sum_ie'_i\ot e''_i)\endalign$$
and (a) follows.

From (a) we deduce:

(b) {\it if $(L',S',\ce')$ is effective (see 21.6) then $\z:\Th(\ce',\ce'')@>>>\bbq^*$ 
is constant.}

\proclaim{Lemma 24.14}If $n\in\Th$ then
$\sum_{y\in S'{}^F}\c_{\ce',\e'}(y)\c_{\ce'',\e''}(n\i yn)$ equals
$$\z(n)q^{\dim S'-\dim L'}|L'{}^F|$$ 
if $n\in\Th(\ce',\ce'')$  and $0$, otherwise.
\endproclaim
In the following proof we write $S$ instead of $S'$ and
$\d$ for the connected component of $N_GL'$ that contains $S$.
By the trace formula for Frobenius
maps our sum is an alternating sum of traces of the Frobenius map on $H^i_c(S,\hce)$ 
where $\hce=\ce'\ot\Ad(n\i)^*\ce''$. If $n\n\Th(\ce',\ce'')$ then $\hce$ has no direct 
summand isomorphic to $\bbq$ hence by 23.5(a), $H^i_c(S',\hce)=0$. Thus we may assume 
that $n\in\Th(\ce',\ce'')$. Then we have canonically $\hce=\bbq\op\hce_1$ as in 24.13. 
As in the proof of 23.5(a) we have $H^i_c(S,\hce_1)=0$ for all $i$ hence 
$H^i_c(S,\hce)=H^i_c(S,\bbq)$. By the definition of $\z(n)$ in 24.13 we see that the 
sum in the lemma is equal to $\z(n)\sum_i(-1)^i\tr(F^*,H^i_c(S,\bbq))$. Let 
$f:\tS@>>>S$ be as in 23.5 (for a fixed $y\in S^F$ and for $N_GL',\d$ instead of 
$G,D$). Then $\tS,f$ are defined over $\FF_q$ and from the proof of 23.5 we see that 
$\tr(F^*,H^i_c(S,\bbq))=\tr(F^*,H^i_c(\tS,\bbq))$. Hence our sum is equal to 
$$\align&\z(n)\sum_i(-1)^i\tr(F^*,H^i_c(\tS,\bbq))=\z(n)|\tS^F|\\&
=\z(n)|{}^\d\cz_{L'}^{0F}||L'{}^F||Z_{L'}(y)^{0F}|\i.\endalign$$
By 10.2, $Z_{L'}(y)^0/{}^\d\cz_{L'}^0$ is a (connected) unipotent group. Hence
$$|{}^\d\cz_{L'}^{0F}||Z_{L'}(y)^{0F}|\i=q^{\dim{}^\d\cz_{L'}^0-\dim Z_{L'}(y)^0}=
q^{\dim S-\dim L'}.$$
The lemma is proved.

\proclaim{Proposition 24.15} In the setup of 24.3 and assuming that $\ce',\ce''$ are 
irreducible and $(L',S',\ce'),(L'',S'',\ce'')$ are effective, the sum 
$$|G^{0F}|\i\sum_{g\in G^F}\c_{\fK',\ph'}(g)\c_{\fK'',\ph''}(g)$$
is $0$ if $\Th(\ce',\ce'')=\em$ and is equal to 
$$\z q^{\dim S'-\dim L'}|\Th(\ce',\ce'')|/|L'{}^F|$$
where $\z=\z(n)$ (see 24.13) for any $n\in\Th(\ce',\ce'')$ if $\Th(\ce',\ce'')\ne\em$.
\endproclaim
This follows from the results in 24.4, 24.13, 24.14.

\proclaim{Proposition 24.16} In the setup of 24.5 let $\ce'=\bbq\bxt\cf'\in\cs(S')$,
$\ce''=\bbq\bxt\cf''\in\cs(S'')$ and let $\e'=1\bxt\bae':F^*\ce'@>\si>>\ce'$,
$\e''=1\bxt\bae'':F^*\ce''@>\si>>\ce''$. Assume that $\cf',\cf''$ are irreducible. Then
the sum 
$$|G^{0F}|\i\sum_{u\in G^F_{un}}Q_{L',G,\boc',\cf',\bae'}(u)
Q_{L'',G,\boc'',\cf'',\bae''}(u)$$
is $0$ if $\Th(\ce',\ce'')=\em$ and is equal to 
$$\z q^{\dim S-\dim L'}|\Th(\ce',\ce'')||L'{}^F|\i|{}^{\d'}\cz_{L'}^{0F}|\i$$
where $\z=\z(n)$ (defined as in 24.13 in terms of $\e',\e''$) for any 
$n\in\Th(\ce',\ce'')$ if $\Th(\ce',\ce'')\ne\em$.
\endproclaim
If $n\in\Th$ (see 24.3) we have clearly
$$\align&\sum_{y\in\boc'{}^F}\c_{\cf',\bae'}(y)\c_{\cf'',\bae''}(n\i yn)\\&
=|{}^{\d'}\cz_{L'}^{0F}|\i\sum_{y\in S'{}^F}\c_{\ce',\e'}(y)\c_{\ce'',\e''}(n\i yn).
\endalign$$
The last sum can be evaluated using the results in 24.13, 24.14. (In our case,
$(L',S',\ce'),(L'',S',\ce'')$ are automatically effective, see 21.8.) We introduce this
in Lemma 24.6. The proposition follows.

\subhead 24.17\endsubhead
Let $A\in\fA(G)$ (see 21.18) and let $\a:F^*A@>\si>>A$ be an isomorphism. We have 
$A=IC(\bY,\ca)$ (extended by $0$ on $G-\bY$) where $\bY$ is the closure of a stratum 
$Y_{L,S}$, $(L,S)\in\AA$ and $\ca$ is an irreducible local system on $\bY$ which is a 
direct summand of $\p_!\tce$ (here $\ce\in\cs(S)$ is irreducible cuspidal and $\p,\tce$
are as in 5.6). By 21.19, 21.20 we can assume that $F(L)=L,F(S)=S$ and 
$F^*\ce\cong\ce$. We consider also $A'\in\fA(G)$ and $\a':F^*A'@>\si>>A'$. Let 
$L',S',\ce',\ca',\bY',\p'_!\tce'$ play the same role for $A'$ as
$L,S,\ce,\ca,\bY,\p_!\tce$ for $A$. (In particular, $A'=IC(\bY',\ca')$ extended by $0$
on $G-\bY'$, $F(L')=L',F(S')=S',F^*\ce'\cong\ce'$.) We fix 
$\e:F^*\ce@>\si>>\ce,\e':F^*\ce'@>\si>>\ce'$. Let 
$$\tG,\G,\fK^w,\ph^w,n_w,g_w,L^w,S^w,\ce^w,\e^w,r,\EE,b_w,\io_i,V_i,\fK_i,\ph_i$$
be associated with $L,S,\ce,\e$ as in 21.6 and let
$$\tG',\G',\fK'{}^{w'},\ph'{}^{w'},n'_{w'},g'_{w'},L'{}^{w'},S'{}^{w'},\ce'{}^{w'},
\e'{}^{w'},r',\EE',b'_{w'},\io'_{i'},V'_{i'},\fK'_{i'},\ph'_{i'}$$
be associated in an analogous way to $L',S',\ce',\e'$. From 21.6(c) we have

(a) $\c_{\fK^w,\ph^w}=\sum_{i\in[1,r]}\tr(b_w\io_i,V_i)\c_{\fK_i,\ph_i}$
\nl
for any effective $w\in\G$. We multiply both sides of (a) by 
$|\G|\i\tr(\io_j\i b_w\i,V_j)$ and sum over all effective $w\in\G$. Using 20.4(c) we 
obtain
$$\c_{\fK_j,\ph_j}=|\G|\i\sum_{w\in\G;\text{eff.}}
\tr(\io_j\i b_w\i,V_j)\c_{\fK^w,\ph^w}$$
for any $j\in[1,r]$. Similarly,
$$\c_{\fK'_{j'},\ph'_{j'}}=|\G'|\i\sum_{w'\in\G';\text{eff.}}
\tr(\io'_{j'}{}\i b'_{w'}{}\i,V'_{j'})\c_{\fK'{}^{w'},\ph'{}^{w'}}$$
for any $j'\in[1,r']$. It follows that
$$\align&|G^{0F}|\i\sum_{g\in G^F}\c_{\fK_j,\ph_j}(g)\c_{\fK'_{j'},\ph'_{j'}}(g)\\&
=|\G|\i|\G'|\i\sum_{w\in\G;\text{eff.}}\sum_{w'\in\G';\text{eff.}}
\tr(\io_j\i b_w\i,V_j)\tr(\io'_{j'}{}\i b'_{w'}{}\i,V'_{j'})\\&\T
|G^{0F}|\i\sum_{g\in G^F}\c_{\fK^w,\ph^w}(g)\c_{\fK'{}^{w'},\ph'{}^{w'}}(g)\tag b
\endalign$$
for any $j\in[1,r],j'\in[1,r']$.

\proclaim{Proposition 24.18} Assume that $\che\ce,\ce'$ are strongly cuspidal. Assume 
also that either $L,L'$ are not $G^0$-conjugate or $\che\ce,\ce'$ are clean. Then
$$|G^{0F}|\i\sum_{g\in G^F}\c_{A,\a}(g)\c_{A',\a'}(g)$$
is $0$ if $A'\not\cong\che A$ and is $q^{\dim S-\dim L}$ if $A'=\che A$ and 
$\a'=\che\a$.
\endproclaim
We may assume that $A=\fK_j,\a=\ph_j,A'=\fK'_{j'},\a'=\ph'_{j'}$ for some $j\in[1,r]$,
$j'\in[1,r']$. If $(L,S,\che\ce),(L',S',\ce')$ are not $G^0$-conjugate, the result
follows from 24.17(b) since by 24.15 we have 
$$|G^{0F}|\i\sum_{g\in G^F}\c_{\fK^w,\ph^w}(g)\c_{\fK'{}^{w'},\ph'{}^{w'}}(g)=0$$ 
for any effective $w,w'$. 

Thus we may assume that $(L,S,\che\ce),(L',S',\ce')$ are $G^0$-conjugate. We may also 
assume that $L=L',S=S',\che\ce=\ce',\e'=\che\e$. We have $\tG'=\tG,\G'=\G$,
$\tce'=(\tce)\che{}$, $\p'_!\tce'=(\p_!\tce)\che{}$. Taking transpose we may identify 
$\EE'=\EE$ as vector spaces but with opposed algebra structures. Then 
$\EE'_w=\EE_{w\i}$ for any $w\in\G$ and we may assume that $b'_w\in\EE'_w$ corresponds 
to $b_w\i\in\EE_{w\i}$. The simple $\EE'$ modules are $V'_i=V_i^*$ (dual space) with 
$\EE'$-action given by taking transpose. We take $\io'_i$ to be the transpose inverse 
of $\io_i$. We have $r'=r,\fK'_i=(\fK_i)\che{},\ph'_i=(\ph_i)\che{}$. We take 
$n'_w=n_w,g'_w=g_w$. Then $L'{}^w=L^w,S'{}^w=S^w,\ce'{}^w=(\ce^w)\che{}$,
$\e'{}^w=(\e^w)\che{}$. Moreover for $w\in\G$, $(L^w,S^w,\ce^w)$ is effective if and 
only if $(L^w,S^w,(\ce^w)\che{})$ is effective.

Let $w,w'\in\G$. We compute 
$$a_{w,w'}=\sh(n\in G^{0F};n\i L^wn=L^{w'},n\i S^wn=S^{w'},
\Ad(n\i)^*(\ce^{w'})\che{}\cong(\ce^w)\che{}).$$
Setting $g_w\i ng_{w'}=\nn$ we see that
$$\align&a_{w,w'}=\sh(\nn\in\tG;F(g_w\nn g_{w'}\i)=g_w\nn g_{w'}\i)\\&=\sh(\nn\in\tG;
n_w\i F(\nn)n_{w'}=\nn)=\sum\Sb y\in\G\\w\i F(y)w'=y\endSb|\D_y^{F_0}|\endalign$$
where $\D_y$ is the $L$-coset in $\tG$ defined by $y$ and $F_0:\D_y@>>>\D_y$ is 
$F_0(\nn)=n_w\i F(\nn)n_{w'}$. Now $L^w$ acts freely on $\D_y$ by 
$l_1:\nn\m l_1*\nn=g_w\i l_1g_w\nn$; this action satisfies 
$F_0(l_1*\nn)=F(l_1)*F_0(\nn)$. It follows that $|\D_y^{F_0}|=|L^{wF}|$ and
$a_{w,w'}|L^{wF}|\i=\sh(y\in\G;w\i F(y)w'=y)$. Using this and 24.15 we see that for any
effective $w,w'\in\G$:
$$\align&|G^{0F}|\i\sum_{g\in G^F}\c_{\fK^w,\ph^w}(g)\c_{\fK'{}^{w'},\ph'{}^{w'}}(g)\\&
=q^{\dim S-\dim L}a_{w,w'}|L^{wF}|\i
=\z_{w,w'}q^{\dim S-\dim L}\sh(y\in\G;w\i F(y)w'=y).\endalign$$
Here $\z_{w,w'}\in\bbq^*$ is defined as follows (assuming that 
$\{y\in\G;w\i F(y)w'=y\}\ne\em$): for any $n\in G^{0F}$ such that
$$n\i L^wn=L^{w'},n\i S^wn=S^{w'},\Ad(n\i)^*(\ce^{w'})\che{}\cong(\ce^w)\che{},$$
$\e^w,(\e^{w'})\che{}$ induce an isomorphism 

$F^*(\ce^w\ot\Ad(n\i)^*(\ce^{w'})\che{})@>>>\ce^w\ot\Ad(n\i)^*(\ce^{w'})\che{}$
\nl
which on the summand $\bbq$ of $\ce^w\ot\Ad(n\i)^*(\ce^{w'})\che{}$ is $\z_{w,w'}$ 
times 
the obvious isomorphism $F^*\bbq@>\si>>\bbq$. 

For $n$ as above we set $g_w\i ng_{w'}=\nn\in\tG$ and let $y$ be the image of $\nn$ in 
$\G$. Let $g\in S$, let $(e_h)$ be a basis of $\ce_g$ and let $(\che e_h)$ be the dual 
basis of $\che\ce_g$. We have an isomorphism $b_y:\Ad(\nn\i)^*\ce@>\si>>\ce$. Taking 
transpose inverse we obtain an isomorphism 
${}^tb_y\i:\Ad(\nn\i)^*\che\ce@>\si>>\che\ce$. This restricts to 
${}^tb_y\i:\che\ce_{\nn\i g\nn}@>\si>>\che\ce_g$ and $\sum_he_h\ot{}^tb_y(\che e_h)$
belongs to the fibre of the summand $\bbq$ of $\ce\ot\Ad(\nn\i)^*\che\ce$ at $g$. We
have
$$\ce^w\ot\Ad(g_{w'}\nn\i g_w\i)^*(\ce^{w'})\che{}=
\Ad(g_w\i)^*(\ce\ot\Ad(\nn\i)^*\che\ce)$$
hence $\sum_he_h\ot{}^tb_y(\che e_h)$ belongs to the fibre of the summand $\bbq$ of 
$$\ce^w\ot\Ad(g_{w'}\nn\i g_w\i)^*(\ce^{w'})\che{}$$
at $g_wgg_w\i\in S^w$. Assuming that $g_wgg_w\i\in S^{wF}$ we see from the definitions 
that
$$\sum_h\e^w(e_h)\ot\che\e^{w'}({}^tb_y(\che e_h))=
\z_{w,w'}\sum_he_h\ot{}^tb_y(\che e_h)\in\ce_g\ot\che\ce_{\nn\i g\nn}$$
hence
$$\sum_h\e(b_we_h)\ot\che\e({}^tb_{w'}\i{}^tb_y(\che e_h))
=\z_{w,w'}\sum_he_h\ot{}^tb_y(\che e_h)\in\ce_g\ot\che\ce_{\nn\i g\nn}.$$
Applying $\e\i\ot\che\e\i$ to both sides and using
$\che\e\i{}^tb_y={}^t(\io\i(b_y))\che\e\i$ we obtain
$$\sum_h b_we_h\ot{}^tb_{w'}\i{}^tb_y(\che e_h)
=\z_{w,w'}\sum_h\e\i(e_h)\ot{}^t(\io\i(b_y))\che\e\i(\che e_h).$$
Setting $e'_h=\e\i(e_h),\che e'_h=\che\e\i(\che e_h)$ we see that $(e'_h),(\che e'_h)$ 
are dual bases of $\ce_{F(g)},\che\ce_{F(g)}$ and we have
$$\sum_h b_we_h\ot{}^tb_{w'}\i{}^tb_y(\che e_h)=\z_{w,w'}\sum_he'_h\ot{}^t(\io\i(b_y))
\che e'_h\in\ce_{F(g)}\ot\che\ce_{F(\nn\i g\nn)}.$$
Applying $1\ot{}^t(\io\i(b_y))\i$ to both sides gives
$$\sum_h b_we_h\ot{}^t(\io\i(b_y))\i{}^tb_{w'}\i{}^tb_y(\che e_h))
=\z_{w,w'}\sum_he'_h\ot\che e'_h\in\ce_{F(g)}\ot\che\ce_{F(g)}.$$
Let $e''_h=b_we_h,\che e''_h={}^tb_w\i\che e_h$. Then $(e''_h),(\che e''_h)$ are dual 
bases of 
$$\ce_{n_wgn_w\i}=\ce_{F(g)},\che\ce_{n_wgn_w\i}=\che\ce_{F(g)}$$
hence $\sum_he'_h\ot\che e'_h=\sum_he''_h\ot\che e''_h$. We see that
$$\sum_h b_we_h\ot{}^t(\io\i(b_y))\i{}^tb_{w'}\i{}^tb_y(\che e_h)
=\z_{w,w'}\sum_h b_we_h\ot {}^tb_w\i\che e_h.$$
This shows that 
$${}^t(\io\i(b_y))\i{}^tb_{w'}\i{}^tb_y(\che e_h)=\z_{w,w'}{}^tb_w\i\che e_h
\in\che\ce_{n_wgn_w\i}$$
for all $h$. Thus ${}^t(\io\i(b_y))\i{}^tb_{w'}\i{}^tb_y=\z_{w,w'}{}^tb_w\i$ as linear
maps $\che\ce_{n_wgn_w\i}@>>>\che\ce_g$ hence $\io\i(b_y)b_{w'}b_y\i=\z_{w,w'}\i b_w$
as linear maps $\ce_{n_wgn_w\i}@>>>\ce_g$. It follows that 
$$\io\i(b_y)b_{w'}b_y\i=\z_{w,w'}\i b_w\in\EE$$
for any $y\in\G$ such that $w\i F(y)w'=y$. In our case, 24.17(b) becomes
$$\align&|G^{0F}|\i\sum_{g\in G^F}\c_{\fK_j,\ph_j}(g)\c_{\fK'_{j'},\ph'_{j'}}(g)\\&
=|\G|^{-2}\sum_{w,w'\in\G;\text{eff.}}\tr(b_w\i\io_j\i,V_j)\\&
\T\tr(\z_{w,w'}b_{w'}\io_{j'},V_{j'})q^{\dim S-\dim L}\sh(y\in\G;w\i F(y)w'=y)
\\&=q^{\dim S-\dim L}|\G|^{-2}\sum_{w\in\G;\text{eff.}}\sum_{y\in\G}
\tr(b_w\i\io_j\i,V_j)\tr(\io\i(b_y\i)b_wb_y\io_{j'},V_{j'}).\endalign$$
We have 
$$\align&\tr(\io\i(b_y\i)b_wb_y\io_{j'},V_{j'})=\tr(\io_{j'}\io\i(b_y\i)b_wb_y,V_{j'}),
\\&\tr(b_y\i\io_{j'}b_wb_y,V_{j'})=\tr(\io_{j'}b_w,V_{j'}),\endalign$$
hence
$$\align&|G^{0F}|\i\sum_{g\in G^F}\c_{\fK_j,\ph_j}(g)\c_{\fK'_{j'},\ph'_{j'}}(g)\\&
=q^{\dim S-\dim L}|\G|\i\sum_{w\in\G;\text{eff.}}
\tr(b_w\i\io_j\i,V_j)\tr(\io_{j'}b_w,V_{j'})=q^{\dim S-\dim L}\d_{j,j'}\endalign$$
where the last step follows from 20.4(c),(b). The proposition is proved.

\subhead 24.19\endsubhead
Assume that we are in the setup of 24.17. Let $\d$ (rep. $\d'$) be the connected 
component of $N_GL$ (resp. $N_GL'$) that contains $S$ (resp. $S'$). Assume that 
$S={}^\d\cz_L^0\T\boc,S'={}^{\d'}\cz_{L'}^0\T\boc'$ where $\boc$ is a unipotent 
$L$-conjugacy class and $\boc'$ is a unipotent $L'$-conjugacy class. Assume that 
$\ce=\bbq\bxt\cf,\ce'=\bbq\bxt\cf'$ where $\cf$ (resp. $\cf'$) is an irrreducible 
$L$- (resp. $L'$-)equivariant local system on $\boc$ (resp. $\boc'$). Let 
$\bae:F^*\cf@>\si>>\cf$ (resp. $\bae':F^*\cf'@>\si>>\cf'$) be the restriction of $\e$ 
(resp. $\e'$). For $w\in\G$ let $\d^w$ be the connected component of $N_G(L^w)$ that 
contains $S^w$. Let $j\in[1,r],j'\in[1,r']$. We can state the following variant of 
Proposition 24.18.

\proclaim{Proposition 24.20} Assume that $\che\ce,\ce'$ are strongly cuspidal. Assume 
also that either $L,L'$ are not $G^0$-conjugate or $\che\ce,\ce'$ are clean. 

(a) If $(L,\boc,\che\cf),(L',\boc',\cf')$ are not $G^0$-conjugate then

$|G^{0F}|\i\sum_{u\in G^F_{un}}\c_{\fK_j,\ph_j}(u)\c_{\fK'_{j'},\ph'_{j'}}(u)=0$.

(b) Assume that $L=L',\boc=\boc',\che\cf=\cf',\che{\bae}=\bae',
\fK'_{j'}=(\fK_{j'})\che{},\ph'_{j'}=(\ph_{j'})\che{}$. Then 
$$\align&|G^{0F}|\i\sum_{u\in G^F_{un}}\c_{\fK_j,\ph_j}(u)\c_{\fK'_{j'},\ph'_{j'}}(u)
\\&=|\G|\sum_{w\in\G}\tr(\io_j\i b_w\i,V_j)\tr(b_w\io_{j'},V_{j'})
|{}^{\d^w}\cz_{L^w}^{0F}|\i q^{\dim S-\dim L}.\endalign$$
\endproclaim
As in 24.17 we have
$$\align&|G^{0F}|\i\sum_{u\in G^F_{un}}\c_{\fK_j,\ph_j}(u)\c_{\fK'_{j'},\ph'_{j'}}(u)
\\&=|\G|\i|\G'|\i\sum_{w\in\G}\sum_{w'\in\G'}\tr(\io_j\i b_w\i,V_j)
\tr(\io'_{j'}{}\i b'_{w'}{}\i,V'_{j'})\\&\T
|G^{0F}|\i\sum_{u\in G^F_{un}}\c_{\fK^w,\ph^w}(u)\c_{\fK'{}^{w'},\ph'{}^{w'}}(u).
\endalign$$
(All elements of $\G,\G'$ are effective in this case.) In the setup of (a), we have
$$|G^{0F}|\i\sum_{u\in G^F_{un}}\c_{\fK^w,\ph^w}(u)\c_{\fK'{}^{w'},\ph'{}^{w'}}(u)=0$$
by 24.16 and (a) follows.

In the setup of (b) we have $\G=\G'$. As in the proof of 24.18 we have (using 24.16 
instead of 24.15):
$$\align&
|G^{0F}|\i\sum_{u\in G^F_{un}}\c_{\fK^w,\ph^w}(u)\c_{\fK'{}^{w'},\ph'{}^{w'}}(u)\\&
=\z_{w,w'}q^{\dim S-\dim L}|{}^{\d^w}\cz_{L^w}^{0F}|\i\sh(y\in\G;w\i F(y)w'=y).
\endalign$$
Hence
$$\align&|G^{0F}|\i\sum_{u\in G^F_{un}}\c_{\fK_j,\ph_j}(u)\c_{\fK'_{j'},\ph'_{j'}}(u)
=|\G|^{-2}\sum_{w,w'\in\G}\tr(\io_j\i b_w\i,V_j)\tr(\io_{j'}b_{w'},V_{j'})\\&\T
\z_{w,w'}q^{\dim S-\dim L}|({}^{\d^w}\cz_{L^w}^0)^F|\i\sh(y\in\G;w\i F(y)w'=y).
\endalign$$
As in the proof of 24.18 we may replace $\tr(\z_{w,w'}\io_{j'}b'_{w'},V_{j'})$ by 
$\tr(\io_{j'}b_w,V_{j'})$ and (b) follows.

\head 25. Properties of cohomology sheaves\endhead
\subhead 25.1\endsubhead
Let $D$ be a connected component of $G$. Let $(L,S),(L',S')\in\AA$. Assume that 
$S\sub D,S'\sub D$. Let $\p:\tY_{L,S}@>>>Y_{L,S}=Y$ be as in 3.13 and let 
$\p':\tY_{L',S'}@>>>Y_{L',S'}$ be the map analogous to $\p$. For any $\ce\in\cs(S)$, 
the local system $\p_!\tce$ on $Y$ is defined as in 5.6; similarly, for any 
$\ce'\in\cs(S')$, the local syatem $\p'_!\tce'$ on $Y_{L',S'}$ is defined.

\proclaim{Proposition 25.2}Let $\ce'\in\cs(S')$. Let $\fK\in\cd(D)$ be 
$IC(\bY_{L',S'},\tce')$, extended by zero on $D-\bY_{L',S'}$. Let $A\in\cd(D)$ be a 
direct summand of $\fK$. Then for any $i\in\ZZ$, there exists $\ce\in\cs(S)$ such that 
the constructible sheaf $\ch^iA|_Y$ is a local system isomorphic to a direct summand of
the local system $\p_!\tce$.
\endproclaim
The proof is given in 25.9.

\subhead 25.3\endsubhead
Let $\cz={}^D\cz_{G^0}^0$. We show:

(a) {\it Let $\ce$ be a $G^0$-equivariant local system on an isolated stratum $C$ of 
$D$. Assume that there exists an integer $n\ge 1$ invertible in $\kk$ and a $\cz$-orbit
$F$ in $C$ such that $\ce|_F$ is $\cz$-equivariant for the action $z:f\m z^nf$ on $F$.
Then $\ce\in\cs(C)$.}
\nl
We may assume that $G$ is generated by $D$ and that $\ce$ is indecomposable. 
We have a canonical direct sum decomposition of $\ce$ such that each summand restricted
to any $\cz$-orbit in $C$ is an isotypical local system isomorphic to a direct sum of 
copies of a fixed $\cz$-equivariant local system of rank $1$. Since $\ce$ is 
indecomposable, it is equal to one of these summands. Thus, if $F$ is a $\cz$-orbit in 
$C$, we may assume that $\ce|_F=\cl^{\op k}$ where $\cl$ is a $\cz$-equivariant local 
system of rank $1$ on $F$. As in 5.3 we can find $\ce_1\in\cs(C)$ of rank $1$ such that
$\ce_1|_F\cong\cl$. Let $\uce=\ce\ot\ce_1^*$. Then $\uce$ is a $G^0$-equivariant local 
system on $C$ such that $\uce|_F\cong\bbq^{\op k}$. Since $G^0$ permutes transitively 
the fibres of $C@>>>C'$ (as in 5.3) we see that the restriction of $\uce$ to any fibre 
of $C@>>>C'$ is $\cong\bbq^{\op k}$. As in 5.3, we see that there is a well defined 
local system $\ce'$ on $C'$ whose inverse image under $C@>>>C'$ is $\uce$. Moreover, 
$\ce'$ is automatically $G^0$-equivariant. It follows that $\uce\in\cs(C)$. Since 
$\ce=\uce\ot\ce_1$ and $\ce_1\in\cs(C)$ we see that $\ce\in\cs(C)$.

\subhead 25.4\endsubhead
Let $n\ge 1$ be such that $\ce'\in\cs_n(S')$. From the definitions we see that

(a) $\fK$ is $\cz\T G^0$-equivariant for the action $(z,x)\m g\m xz^ngx\i$ on $D$.
\nl
Let $\ch=\ch^i\fK$, a constructible sheaf on $D$. Let $\d$ be the connected component 
of $N_GL$ such that $S\sub\d$. Let $\ct={}^\d\cz_L^0$ .

We have a canonical map $S@>>>S_s$, $g\m g_s$ and $S_s$ is a single orbit for the 
$\ct\T L$-action $(z,x):y\m xzyx\i$. Let $\ss$ be the set of $\ct$-cosets in $S_s$. 
Then $L$ acts transitively on $\ss$. We fix a $\ct$-coset $\t$ on $S_s$. Let 
$R=\{g\in S;g_s\in\t\}$, $R^*=\{g\in S^*;g_s\in\t\}=R\cap S^*$. Then $R^*$ is open in 
$R$. By the proof of 3.11, $R^*$ is dense in $R$. Now $R$ is a single orbit for the 
group $\ct\T N_L\t$ acting by $(z,l):g\m lzgl\i$. Hence $R$ is smooth of pure 
dimension. Let $s\in\t$. Then $Z_L(s)^0$ is independent of the choice of $s$: it is 
$N_L(\t)^0=Z_L(\t)^0$. Hence R is union of finitely many orbits of $\ct\T Z_L(\t)^0$.

\proclaim{Lemma 25.5}Let $s\in\t$. Then 

(a) $R\sub Z_G(s)$.

(b) Any connected component of $s\i R$ contains some unipotent element.
\endproclaim
We can find $y\in S$ with $y_s=s$. We have $y\in R$ and $ys=sy$. Let $y'\in R$. We have
$y'=lzyl\i$ where $l\in L,l\i sl=z's$, $z,z'\in\ct$. We must show that $sy'=y's$ or 
that $slzyl\i=lzyl\i s$ or that $l\i slzy=zyl\i sl$ or that $z'szy=zyz's$ or that
$z'zsy=zz'ys$ or that $sy=ys$. This proves (a). We prove (b). Let $g\in R$. Then 
$g_s=zs$ for some $z\in\ct$. Let 
$g'=z\i g$. Then $g'_s=z\i g_s=s$. Hence $g'=su$ where $u$ is unipotent in $Z_G(s)$. 
Now $\ct R=R$ hence $g',g$ are in the same connected component of $R$. Hence the 
connected component of $g$ in $R$ contains an element $g'$ such that $s\i g'$ is 
unipotent. This proves (b).

\proclaim{Lemma 25.6}There exists a local system $\cf$ on $R$ and an integer $n\ge 1$
invertible in $\kk$ such that $\cf$ is $\ct\T Z_L(\t)^0$-equivariant for the action 
$(z,x)\m g\m xz^ngx\i$ and $\cf|_{R^*}\cong\ch|_{R^*}$.
\endproclaim
Let $R_1$ be a connected component of $R$ which is a $\ct\T Z_L(s)^0$-orbit in $R$. Let
$R_1^*=R^*\cap R$, an open dense subset of $R_1$. It is enough to show: 

{\it there exists a local system $\cf_1$ on $R_1$ and an integer $n\ge 1$ invertible in
$\kk$ such that $\cf_1$ is $\ct\T Z_L(\t)^0$-equivariant for the action 
$(z,x)\m g\m xz^ngx\i$ and $\cf_1|_{R^*_1}\cong\ch|_{R^*_1}$.}
\nl
Let $\t^*$ be the set of all $s\in\t$ such that $s=y_s$ for some $y\in S^*$. By the 
proof of 3.11 we have $\t^*\ne\em$. Let $s\in\t^*$. Let $\d_1$ be the connected 
component of $Z_{N_GL}(s)$ such that $R_1\sub s\d_1$. By 25.5(b), $\d_1$ contains some 
unipotent element. By 16.12(b) we can find an open subset $\cu$ of
$\d_1$ containing all unipotents of $\d_1$, such that $\e^*(\fK|_{s\cu})$ is isomorphic
to the restriction to $\cu$ of $\tK$ where $\tK$ is a direct sum of finitely many 
objects $\tK^k\in\cd(\d_1)$ each one of the same type as $\fK$ (for $Z_G(s)$ instead of
$G$). (Here $\e:\cu@>>>s\cu$ is $g\m sg$.) By an analogue of 25.4(a), each $\tK^k$ is 
${}^{\d_1}\cz_{Z_G(s)^0}^0\T Z_G(s)^0$-equivariant hence $\ct\T Z_L(\t)^0$-equivariant 
for an action as in the Lemma. (We have $Z_G(s)^0=Z_L(s)^0=Z_L(\t)^0$, since $s\in\t^*$
and ${}^{\d_1}\cz_{Z_G(s)^0}^0={}^{\d_1}\cz_{Z_L(s)^0}^0=T_{N_GL}(y)={}^\d\cz_L^0=\ct$,
since $y\in R_1,y_s=s$ is isolated in $N_GL$.) Then $\ch^i\tK$ is a 
$\ct\T Z_L(\t)^0$-equivariant constructible sheaf on $\d_1$ whose restriction to $\cu$ 
is $\cong\ch^i(\e^*\fK|_{s\cu})$. Let $\cf_1$ be the inverse image of 
$\ch^i\tK|_{s\i R_1}$ under $R_1@>>>s\i R_1,g\m s\i g$. This is a 
$\ct\T Z_L(\t)^0$-equivariant constructible sheaf on $R_1$. Since $R_1$ is a single 
orbit, we see that $\cf_1$ is a $\ct\T Z_L(\t)^0$-equivariant local system on $R_1$. We
have $\cf_1|_\cv\cong\ch|_\cv$ where $\cv=R_1\cap s\cu$, an open subset of $R_1$ 
containing $\{g\in R_1;g_s=s\}$. We summarize: for any $s\in\t^*$ there exists an open 
subset $\cv(s)$ of $R_1$ and a $\ct\T Z_L(\t)^0$-equivariant local system $\cf_s$ on 
$R_1$ such that $\{g\in R_1;g_s=s\}\sub\cv(s)$ and $\cf_s|_{\cv(s)}\cong\ch|_{\cv(s)}$.
If $g\in R_1^*$ then $g\in R_1$ and $g_s\in\t^*$ hence $g\in\cv(g_s)$. Thus 
$R_1^*\sub R'$ where $R'=\cup_{s\in\t^*}\cv(s)\sub R_1$. The constructible sheaf 
$\ch|_{R'}$ is a local system when restricted to any of the open sets $\cv(s)$ that 
cover $R'$. Hence $\ch|_{R'}$ is a local system. Let $s\in\t^*$ be such that 
$\cv(s)\ne\em$. The local systems $\ch|_{R'},\cf_s|_{R'}$ have isomorphic restrictions 
to $\cv(s)$ hence $\ch|_{R'}\cong\cf_s|_{R'}$. (Note that $R'$ is smooth, irreducible, 
since $R_1$ is smooth, irreducible. Hence 
$\ch|_{R'}=IC(R',\ch|_{\cv(s)})\cong IC(R',\cf_s|_{\cv(s)})=\cf_s|_{R'}$.) Since 
$R_1^*\sub R'$, it follows that $\ch|_{R_1^*}\cong\cf_s|_{R_1^*}$. We set $\cf=\cf_s$.
Then $\cf$ has the required properties. The lemma is proved.

\proclaim{Lemma 25.7}There exists $\ce\in\cs(S)$ such that $\ce|_{S^*}\cong\ch|_{S^*}$.
\endproclaim
Define $\r:S@>>>\ss$ by $g\m(\ct-\text{coset of }g_s)$. Let $\r':S^*@>>>\ss$ be
the restriction of $\r$. Then $\r,\r'$ are $L$-equivariant maps with $L$ acting 
transitively on $\ss$ and $R$ (resp. $R^*$) is a fibre of $\r$ (resp. $\r'$). Since 
$\ch|_{S^*}$ is an $L$-equivariant constructible sheaf whose restriction to any fibre 
of $\r'$ is a local system (see 25.6), we see that $\ch|_{S^*}$ is an $L$-equivariant 
local system. Hence $\un K:=IC(S,\ch|_{S^*})$ is a well defined object 
of $\cd(S)$. Let $\un K'=IC(R,\ch|_{R^*})$. Using the $L$-homogeneity of $\ss$ we see 
that for any $j\in\ZZ$ we have $\ch^j\un K|_R=\ch^j\un K'$ and $\ch^j\un K$ is an 
$L$-equivariant constructible sheaf. From 25.6 we see that $\un K'$ is a local
system on $R$. In particular, $\ch^j\un K'=0$ for $j>0$. Since for $j>0$, $\ch^j\un K$ 
is an $L$-equivariant constructible sheaf on $S$ whose restriction to $R$ is $0$, we 
see that $\ch^j\un K=0$. Similarly, since $\ch^0\un K$ is an $L$-equivariant 
constructible sheaf on $S$ whose restriction to $R$ is a local system, we see that 
$\ch^0\un K$ is a local system. Thus, $\un K=\ch^0\un K$ is a local system $\ce$ on $S$
whose restriction to $S^*$ is $\ch|_{S^*}$. Since $\ce=IC(S,\ch|_{S^*})$ and 
$\ch|_{S^*}$ is $L$-equivariant, we see that $\ce$ is $L$-equivariant. Since 
$\ce|_R=\un K'$ and $\un K'$ is $\ct$-equivariant, we may apply 25.3(a) with $G,C$ 
replaced by $N_GL,S$ and we see that $\ce\in\cs(S)$. The lemma is proved.

\proclaim{Lemma 25.8}Let $\tY=\tY_{L,S}$. Define $\x:\tY@>>>D$ by $\x(g,xL)=g$. Then 
$\x^*\ch\cong\tce$ where $\tce$ is defined in terms of $\ce\in\cs(S)$ of Lemma 25.7.
\endproclaim
Let $\hY=\{(g,x)\in D\T G^0;x\i gx\in S^*\}$. Define $b:\hY@>>>S$ and $b':\hY@>>>S^*$  
by $(g,x)\m x\i gx$. Define $c:\hY@>>>D$ by $c(g,x)=g$. By the definition of $\ce$ we
have $b^*\ce=b'{}^*(\ce|_{S^*})\cong b'{}^*(\ch|_{S^*})$. Define $r:\hY@>>>G^0\T D$ by 
$(g,x)\m(x,x\i gx)$. Define $p_1,p_2:G^0\T D@>>>D$ by $p_1(x,d)=d,p_2(x,d)=xdx\i$. 
Since $\ch$ is a $G^0$-equivariant constructible sheaf on $D$, we have 
$p_1^*\ch\cong p_2^*\ch$, hence $r^*p_1^*\ch\cong r^*p_2^*\ch$. Thus, 
$b'{}^*(\ch|_{S^*})\cong c^*\ch$. It follows that $c^*\ch\cong b^*\ce$. As in 5.6,
define $a:\hY@>>>\tY$ by $(g,x)\m(g,xL)$. We have $c^*\ch=a^*(\x^*\ch),b^*\ce=a^*\tce$.
Hence $a^*(\x^*\ch)\cong a^*\tce$. Since $a^*\tce$ is a local system, we see that 
$a^*(\x^*\ch)$ is a local system. Since $a$ is a principal $L$-bundle it follows that
$\x^*\ch$ is a local system and that $\x^*\ch\cong\tce$. The lemma is proved.

\subhead 25.9\endsubhead
We now prove Proposition 25.2. We may assume that $A=\fK$. Let $j:Y@>>>D$ be the 
inclusion. With notation in 25.8 we have $\p^*j^*\ch=\x^*\ch$ hence $\p^*j^*\ch=\tce$. 
Since $\p$ is a finite unramified covering and $\tce$ is a local system, it follows 
that $j^*\ch$ is a local system. Since $\bbq$ is a direct summand of $\p_!\p^*\bbq$ we 
see that $j^*\ch$ is a direct summand of 
$(j^*\ch)\ot(\p_!\p^*\bbq)=\p_!\p^*(j^*\ch)=\p_!\tce$. Proposition 25.2 is proved.

\proclaim{Lemma 25.10}Let $(L,S),(L',S')\in\AA$. Assume that $Y_{L',S'}\sub\bY_{L,S}$. 

(a) For any $L$-conjugacy class $C$ in $S$ and any $L'$-conjugacy class $C'$ in $S'$ we
have $\dim L'-\dim C'\ge\dim L-\dim C$.

(b) For any $G^0$-conjugacy class $\boc$ in $Y_{L,S}$ and any $G^0$-conjugacy class 
$\boc'$ in $Y_{L',S'}$ we have $\dim\boc'\le\dim\boc$.
\endproclaim
We prove (a). By 7.2(c) we may assume that $L'=G^0$. Then $S'$ is as isolated stratum 
of $G$ and $S'\sub\bY$ where $Y=Y_{S,L}$. Let $\aa=\aa_{L,S},\s:\bY@>>>\aa$ be as in 
7.2. Set $a=\s(C')\in\aa$. By 7.16(b), $\bY^a$ has pure dimension $\dim G^0/L+\dim C$. 
Since $C'\sub\bY^a$ we have $\dim C'\le\dim G^0/L+\dim C$, as required.

We prove (b). Let $g\in\boc,g'\in\boc'$. We can assume that $g\in S^*,g'\in S'{}^*$.
Let $C$ be the $L$-conjugacy class of $g$; let $C'$ be the $L'$-conjugacy class of 
$g'$. By the definition of $S^*$ we have $Z_G(g)^0=Z_L(g)^0$ hence 
$\dim Z_{G^0}(g)=\dim Z_L(g)$ so that $\dim\boc=\dim G^0/L+\dim C$. Similarly,
$\dim\boc'=\dim G^0/L'+\dim C'$. Hence (b) follows from (a).

\subhead 25.11\endsubhead
For $(L,S)\in\AA$ let $\fT_{L,S}$ be the union of all $G^0$-conjugacy classes in
$\bY_{L,S}$ whose dimension equals the dimension of some/any $G^0$-conjugacy class in
$Y_{L,S}$. (See 3.4.) From 3.15 and 3.4 we see that $\fT_{L,S}$ is a union of strata of
$G$. Clearly,

(a) $Y_{L,S}\sub\fT_{L,S}\sub\bY_{L,S}$.

\proclaim{Lemma 25.12}For $(L,S)\in\AA$, $\fT_{L,S}$ is open dense in $\bY_{L,S}$.
\endproclaim
The fact that it is dense follows from 25.11(a). It remains to show that
$\bY_{L,S}-\fT_{L,S}$ is closed in $\bY_{L,S}$. Since $\bY_{L,S}$ is a union of strata
(see 3.15) and $\fT_{L,S}$ is a union of strata (by definition) it is enough to verify
the following statement:

$Y_{L'',S''}\sub\bY_{L,S}-\fT_{L,S},Y_{L',S'}\sub\bY_{L'',S''}\imp
Y_{L',S'}\sub\bY_{L,S}-\fT_{L,S}$.
\nl
Let $\boc$ (resp. $\boc',\boc''$) be a $G^0$-conjugacy class in $Y_{L,S}$ (resp. 
$Y_{L',S'}$, $Y_{L'',S''}$). Using Lemma 25.10 we see that $\dim\boc''<\dim\boc$,
$\dim\boc'\le\dim\boc''$. Hence $\dim\boc'<\dim\boc$ and 
$Y_{L',S'}\sub\bY_{L,S}-\fT_{L,S}$. The lemma is proved.

\subhead 25.13\endsubhead
Let $(L,S)\in\AA$ and let $P$ be a parabolic of $G^0$ with Levi $L$ such that
$S\sub N_GP$. Let $\fT=\fT_{L,S}$. Let $\ps:X@>>>\bY_{L,S}$ be as in 3.14. For any 
stratum $S'$ of $N_GP\cap N_GL$ such that $S'\sub\bS$ let $X_{S'}$ be as in 5.6. Let 
$\fX=\ps\i(\fT)$ and let $\ps':\fX@>>>\fT$ be the restriction of $\ps$. We show:

(a) $\fX\sub X_S$;

(b) $\ps'$ {\it has finite fibres;}

(c) $\fX$ {\it is smooth.}
\nl
Let $g\in\fT$. We must show that $\ps\i(g)\cap X_{S'}$ is empty if $S'\ne S$ and is
finite if $S'=S$. By 4.4(b), 
$\dim(\ps\i(g)\cap X_{S'})\le(\dim G^0/L-\dim\boc_1+\dim C')/2$ where $\boc_1$ is the 
$G^0$-conjugacy class of $g$ and $C'$ is any $L$-conjugacy class in $S'$. Let $C$ be an
$L$-conjugacy class in $S$. We have $\dim C'\le\dim C$ with strict inequality if 
$S'\ne S$. Hence $\dim(\ps\i(g)\cap X_{S'})\le(\dim G^0/L-\dim\boc_1+\dim C)/2$ with
strict inequality if $S'\ne S$. As in the proof of 25.10 we have 
$\dim G^0/L+\dim C=\dim\boc$ where $\boc$ is a $G^0$-conjugacy class in $Y_{L,S}$. Also
$\dim\boc=\dim\boc_1$ by the choice of $g$. Thus $\dim(\ps\i(g)\cap X_{S'})\le 0$ with 
strict inequality if $S'\ne S$. This completes the proof of (a),(b). 

Using (a) and Lemma 25.12 we see that $\fX$ is an open subset of $X_S$ which is smooth.
Hence (c) holds.

\proclaim{Proposition 25.14}Let $(L,S)\in\AA$, let $\ce\in\cs(S)$ and let 
$\fT=\fT_{L,S}$. Then $IC(\bY_{L,S},\p_!\tce)|_{\fT}$ is a constructible sheaf.
\endproclaim
Let $\bce$ be the local system on $X_S$ defined in 5.6. Using 5.7 it is enough to show
that $\ps_!(IC(X,\bce))|_{\fT}$ is a constructible sheaf or equivalently that
$\ps'_!(IC(X,\bce)|_{\fX})$ is a constructible sheaf. Since $\fX\sub X_S$ (see 
25.13(a)) we have $IC(X,\bce)|_{\fX}=\bce|_{\fX}$ and it is enough to show that
$\ps'_!(\bce|_{\fX})$ is a constructible sheaf. This is clear from 25.13(b).

\head 26. The variety $Z_{J,D}$\endhead
\subhead 26.1\endsubhead
Let $\cb$ be the variety of Borel subgroups of $G^0$. Let $\WW$ be the set of 
$G^0$-orbits on $\cb\T\cb$ ($G^0$ acts by conjugation on both factors). For
$B,B'\in\cb$ we write $\po(B,B')=w$ if the $G^0$-orbit of $(B,B')$ is $w$. There is a
unique group structure on $\WW$ such that whenever $B,B',B''\in\cb$ have a common 
maximal torus, we have $\po(B,B')\po(B',B'')=\po(B,B'')$. Then $\WW$ is a finite 
Coxeter group (called the {\it Weyl group}) with length function $l:\WW@>>>\NN$ which 
attaches to a $G^0$-orbit its dimension minus $\dim\cb$. Let $\le$ be the standard 
partial order of the Coxeter group $\WW$. Let $\II=\{w\in\WW;l(w)=1\}$. For $J\sub\II$ 
let $\WW_J$ be the subset of $\WW$ generated by $J$.

If $P$ is a parabolic of $G^0$, the set of all $w\in\WW$ such that $w=\po(B,B')$ for 
some $B,B'\in\cb,B\sub P,B'\sub P$ is of the form $\WW_J$ for a well defined 
$J\sub\II$; we then say that $P$ has type $J$. For $J\sub\II$ let $\cp_J$ be the set of
all parabolics of type $J$ of $G^0$. For $P\in\cp_J,Q\in\cp_K$ there is a well defined 
element $u=\po(P,Q)\in\WW$ such that $u\le\po(B,B')$ for any $B,B'\in\cb,B\sub P$,
$B'\sub Q$ and $u=\po(B_1,B'_1)$ for some $B_1,B'_1\in\cb,B_1\sub P,B'_1\sub Q$.

\subhead 26.2\endsubhead
In the remainder of this section we fix a connected component $D$ of $G$.

There is a unique isomorphism $\e:\WW@>>>\WW$ such that $\e(\II)=\II$ and such that

$g\in D,P\in\cp_J\imp gPg\in\cp_{\e(J)}$. 
\nl
We fix $J\sub\II$. Following \cite{\BE}, let $\ct(J,\e)$ be the set of all sequences 
$(J_n,w_n)_{n\ge 0}$ where $J_n\sub\II$ and $w_n\in\WW$ are such that
$$J=J_0\supset J_1\supset J_2\supset\do,$$
$$J_n=J_{n-1}\cap\e\i(w_{n-1}J_{n-1}w_{n-1}\i)\text{ for }n\ge 1,$$
$$l(w_ny)>l(w_n)\text{ for all }y\in J_n,n\ge 0,$$
$$l(y'w_n)>l(w_n)\text{ for all }y'\in\e(J_n),n\ge 0,$$
$$w_n\in W_{\e(J_n)}w_{n-1}W_{J_{n-1}}\text{ for }n\ge 1.$$
Then $\ct(J,\e)$ is a finite set.

For $(P,P')\in\cp_J\T\cp_{\e(J)}$ let $A(P,P')=\{g\in D;gPg\i=P'\}$. Let $Z_{J,D}$ be 
the set of all triples $(P,P',\g)$ where $P\in\cp_J,P'\in\cp_{\e(J)}$ and 
$\g\in U_{P'}\bsl A(P,P')=A(P,P')/U_P$. Following \cite{\PCS, 3.11}, to any 
$(P,P',\g)\in Z_{J,D}$ we associate an element $(J_n,w_n)_{n\ge 0}\in\ct(J,\e)$ and 
two sequences of parabolics $P^n,P'{}^n,(n\ge 0)$ by the requirements:
$$P'{}^0=P',P^0=P, P'{}^n=(P^{n-1}\cap P'{}^{n-1})U_{P'{}^{n-1}}\in\cp_{\e(J_n)},$$ 
$$P^n=g\i P'{}^ng\in\cp_{J_n}, g\in\g, w_n=\po(P'_n,P^n).$$
\nl
We write $(J_n,w_n)_{n\ge 0}=\b'(P,P',\g)$. For $\tt\in\ct(J,\e)$ let
$${}^\tt Z_{J,D}=\{(P,P',\g)\in Z_{J,D};\b'(P,P',\g)=\tt\}.$$
Then $({}^\tt Z_{J,D})_{\tt\in\ct(J,\e)}$ is a partition of $Z_{J,D}$ into locally 
closed subvarieties. Now $G^0$ acts on $Z_{J,D}$ by 
$h:(P,P',\g)\m(hPh\i,hP'h\i,h\g h\i)$. This action preserves each of the pieces 
${}^\tt Z_{J,D}$.

\subhead 26.3\endsubhead
Let $\tt=(J_n,w_n)_{n\ge 0}\in\ct(J,\e)$. For $r\gg 0$, $J_r,w_r$ are independent of
$r$; we denote them by $J_\iy,w$. Then 
$$wJ_\iy w\i=\e(J_\iy),l(wy)>l(w),l(y'w)>l(w)\text{ for all }y\in J_\iy,
y'\in\e(J_\iy).$$
Let $R_\tt=\{(\tQ,\tQ',\g')\in Z_{J_\iy,D},\po(\tQ',\tQ)=w\}$. We choose 
$Q\in\cp_{J_\iy},Q'\in\cp_{\e(J_\iy)}$ such that $\po(Q',Q)=w$. We can find a common 
Levi $L$ for $Q$ and $Q'$. Let 
$$C=\{g\in D;gLg\i=L,gQg\i=Q'\}=\{g\in D;gLg\i=L,\po(gQg\i,Q)=w\}.$$
Let $A$ be a simple perverse sheaf on $C$ which
is admissible in the sense of 6.7 (this concept is well defined since $C$ is a 
connected component of the reductive group $N_GL$.) Then $A$ is $L$-equivariant for the
conjugation action of $L$ hence it is also $(Q\cap Q')$-equivariant where $Q\cap Q'$ 
acts via its quotient $(Q\cap Q')/U_{Q\cap Q'}=L$. Hence there is a well defined simple
perverse sheaf $A'$ on $G^0\T_{Q\cap Q'}C$ (here $Q\cap Q'$ acts on $G^0$ by right 
translation) such that $\tj A'=\wt{pr}_2A$ in the obvious diagram
$G^0\T_{Q\cap Q'}C@<j<<G^0\T C@>pr_2>>C$. We may regard $A'$ as a simple perverse sheaf
on $R_\tt$ via the isomorphism
$$G^0\T_{Q\cap Q'}C@>\si>>R_\tt,(g,c)\m(gQg\i,gQ'g\i,gU_{Q'}cU_Qg\i).\tag a$$
Define $\vt_\tt:{}^\tt Z_{J,D}@>>>R_\tt$ by $(P,P',\g)\m(P^r,P'{}^r,\g U_{P^r})$ where 
$r\gg 0$ and $P^r,P'{}^r$ are attached to $(P,P',\g)$ as in 26.2. Now $G^0$ acts on 
$R_\tt$ by $h:(\tQ,\tQ',\g')\m(h\tQ h\i,h\tQ' h\i,h\g'h\i)$ and $\vt_\tt$ is 
$G^0$-equivariant. By \cite{\PCS, 3.12}, $\vt_\tt$ is an iterated affine space bundle. 
Let $\tA=\ti\vt(A')$, a simple perverse sheaf on ${}^\tt Z_{J,D}$. Let $\hA$ be the 
simple perverse sheaf on $Z_{J,D}$ whose support is the closure in $Z_{J,D}$ of 
$\supp(\tA)$ and whose restriction to ${}^\tt Z_{J,D}$ is $\tA$. A simple perverse 
sheaf on ${}^\tt Z_{J,D}$ is said to be {\it admissible} if it is of the form $\tA$ for
some $A$ as above. This concept does not depend on the choice of $Q,Q',L$ since any two
such triples are $G^0$-conjugate. Note that $A\m\tA$ is a bijection between the set of 
isomorphism classes of simple perverse sheaves on $C$ that are admissible and the set 
of isomorphism classes of simple perverse sheaves on ${}^\tt Z_{J,D}$ that are 
admissible. A simple perverse sheaf on $Z_{J,D}$ is said to be {\it admissible} if it 
is of the form $\hA$ for some $\tt\in\ct(J,\e)$ and some $A$ as above. Note that 
$(\tt,A)\m\hA$ is a bijection between the set of pairs consisting of an element of 
$\ct(J,\e)$ and an isomorphism class of a simple perverse sheaf on ${}^\tt Z_{J,D}$ 
that is admissible and the set of isomorphism classes of simple perverse sheaves on 
$Z_{J,D}$ that are admissible. 

When $J=\II$ then $\ct(J,\e)$ consists of a single element $\tt=(J_n,w_n)_{n\ge 0}$ 
where $J_n=\II,w_n=1$ for all $n$. Now $Z_{\II,D}$ consists of all triples 
$(G^0,G^0,g)$ with $g\in D$. We identify $Z_{\II,D}=D$ in the obvious way. A simple
perverse sheaf on $Z_{\II,D}$ is admissible in the sense just defined if and only if it
is admissible on $D$ in the sense of 6.7.

\subhead 26.4\endsubhead
In the remainder of this section we assume that $\kk$ is an algebraic closure of a
finite field $\FF_q$ and that $G$ has a fixed $\FF_q$-rational structure with Frobenius
map $F:G@>>>G$. There are induced maps $F:\cb@>>>\cb,F:\WW@>>>\WW$; the last map
restricts to a bijection $F:\II@>>>\II$. Let $J\sub\II$ be such that $F(J)=J$. Then
$F:G@>>>G$ induces a map $F:\cp_J@>>>\cp_J$. We assume that $F(D)=D$. Then 
$\e:\II@>>>\II$ commutes with $F$ hence $F(\e(J))=\e(J)$. Hence $F:G@>>>G$ induces a 
map $F:\cp_{\e(J)}@>>>\cp_{\e(J)}$. 

For $(P,P')\in\cp_J\T\cp_{\e(J)}$ we have $g\in A(P,P')\imp F(g)\in A(F(P),F(P'))$ and 
$g\m F(g)$ induces a map $F:A(P,P')/U_P@>>>A(F(P),F(P'))/U_{F(P)}$. Define 
$F:Z_{J,D}@>>>Z_{J,D}$ by $F(P,P',\g)=(F(P),F(P'),F(\g))$; this is the Frobenius map 
for an $\FF_q$-rational structure on $Z_{J,D}$. The $G^0$-action on $Z_{J,D}$ restricts
to a $G^{0F}$-action on $Z_{J,D}^F$. Let $\UU$ be the vector space of functions 
$Z_{J,D}^F@>>>\bbq$ that are constant on $G^{0F}$-orbits.

\proclaim{Theorem 26.5}Let $\ca_J$ be a set of representatives for the isomorphism 
classes of admissible simple perverse sheaves $K$ on $Z_{J,D}$ such that $F^*K\cong K$.
For each $K\in\ca_J$ we choose an isomorphism $\a:F^*K@>\si>>K$. The characteristic 
functions $\c_{K,\a}$ (one for each $K\in\ca_J$) form a $\bbq$-basis of $\UU$.
\endproclaim
In this proof we write $Z,{}^\tt Z$ instead of $Z_{J,D},{}^\tt Z_{J,D}$. Define 
$F:\ct(J,\e)@>\si>>\ct(J,\e)$ by $(J_n,w_n)_{n\ge 0}\m(F(J_n),F(w_n))_{n\ge 0}$. For 
any $\tt\in\ct(J,\e)$ we have $F({}^\tt Z)={}^{F(\tt)}Z$. In particular, we have 
$Z^F=\sqc_{\tt\in\ct(J,\e);F(\tt)=\tt}{}^\tt Z^F$
where ${}^\tt Z^F={}^\tt Z\cap Z^F$. It follows that 
$\UU=\op_{\tt\in\ct(J,\e);F(\tt)=\tt}{}^\tt\UU$ where ${}^\tt\UU$ is the vector space 
of functions ${}^\tt Z^F@>>>\bbq$ that are constant on $G^{0F}$-orbits. (We identify 
any such function with a function $Z^F@>>>\bbq$ which is zero on the complement of 
${}^\tt Z^F$.) 

For any integer $t$ let $Z_{\le t}=\cup_{\tt\in\ct(J,\e);|\tt|\le t}{}^\tt Z$ where 
$|\tt|=\dim{}^\tt Z$. This is a closed subvariety of $Z$ since $\cup_\tt{}^\tt Z$ is a 
partition of $Z$ into finitely many locally closed subvarieties. Moreover, if 
$|\tt|=t$ then ${}^\tt Z\cup Z_{\le t-1}$ is a closed subvariety of $Z$.

Let $\UU_{\le t}$ be the vector space of functions $Z_{\le t}^F@>>>\bbq$ that are 
constant on $G^{0F}$-orbits. (We identify any such function with a function 
$Z^F@>>>\bbq$ which is zero on the complement of ${}^\tt Z^F$.) We have 
$\UU_{\le 0}\sub\UU_{\le 1}\sub\do$ and
$\UU_{\le t}=\op_{\tt\in\ct(J,\e);F(\tt)=\tt,|\tt|\le t}{}^\tt\UU$.

Now let $K\in\ca_J$. Since the sets ${}^\tt Z\cap\supp(K)$ form a partition of 
$\supp(K)$ into finitely many locally closed subsets, there is a unique
$\tt\in\ct(J,\e)$ such that ${}^\tt Z\cap\supp(K)$ is open dense in $\supp(K)$. Since 
$F^*K\cong K$ we have necessarily $F(\tt)=\tt$. Let $t=|\tt|$. Since 
${}^\tt Z\cap\supp(K)\sub{}^\tt Z$ and ${}^\tt Z\cup Z_{\le t-1}$ is closed, we see 
that the closure of ${}^\tt Z\cap\supp(K)$ is contained in ${}^\tt Z\cup Z_{\le t-1}$. 
Thus, $\supp(K)\sub{}^\tt Z\cup Z_{\le t-1}$. Hence we can write uniquely
$\c_{K,\a}=\c'_{K,\a}+\c''_{K,\a}$ where $\c'_{K,\a}$ is the restriction of $\c_{K,\a}$
to ${}^\tt Z$ (extended by $0$ on $Z-{}^\tt Z$) and $\c''_{K,\a}\in\UU_{\le t-1}$. Note
that $\c_{K,\a},\c'_{K,\a}$ are contained in $\UU_{\le t}$ and they are equal modulo 
$\UU_{\le t-1}$. Thus, in order to prove that the functions $\c_{K,\a}$ form a basis of
$\UU$ it is enough to prove that the functions $\c'_{K,\a}$ form a basis of $\UU$. More
precisely, we will show that, for any $\tt\in\ct(J,\e)^F$, the functions $\c'_{K,\a}$ 
with $K\in\ca_J$ such that ${}^\tt Z\cap\supp(K)$ is open dense in $\supp(K)$, form a
basis of ${}^\tt\UU$. An equivalent statement is:

(a) {\it For any $\tt\in\ct(J,\e)^F$ let ${}^\tt\ca_J$ be a set of representatives for 
the isomorphism classes of admissible simple perverse sheaves $K'$ on ${}^\tt Z$ such 
that $F^*K\cong K$. For each $K'\in{}^\tt\ca_J$ we choose an isomorphism 
$\a':F^*K'@>\si>>K'$. Then the characteristic functions $\c_{K',\a'}$ (one for each 
$K'\in{}^\tt\ca_J$) form a $\bbq$-basis of ${}^\tt\UU$.}
\nl
Let $R=R_\tt,Q,Q',L,C$ be as in 26.3. Then $R$ is defined over $\FF_q$, with 
Frobenius map $F:R_\tt@>>>R_\tt,(\tQ,\tQ',\g')\m(F(\tQ),F(\tQ'),F(\g'))$. We may assume
that $F(Q)=Q,F(Q')=Q',F(L)=L'$. Then $F(C)=C$. The map $\vt_\tt:{}^\tt Z@>>>R$ in 26.3
is $G^0$-equivariant and commutes with $F$. 

We show that $\vt_\tt$ induces a bijection between the set of $G^{0F}$-orbits on 
${}^\tt Z$ and the set of $G^{0F}$-orbits on $R$. Now $\vt_\tt$ is a composition 
$${}^\tt Z_{J,D}@>{}^1\vt>>{}^{\tt_1}Z_{J_1,D}@>{}^2\vt>>\do@>{}^r\vt>>
{}^{\tt_r}Z_{J_r,D}=R$$
where $\tt_i=(J_n,w_n)_{n\ge i}$, ${}^1\vt(P,P',\g)=(P^1,P'{}^1,\g U_{P^1})$ (notation 
of 26.2) and ${}^i\vt$ is defined for $i\ge 2$ just like ${}^1\vt$ with $\tt_{i-1}$ 
instead of $\tt$. Each ${}^i\vt$ is $G^0$-equivariant and commutes with $F$. Hence it
is enough to show that ${}^i\vt$ induces a bijection between the set of $G^{0F}$-orbits
on ${}^{\tt_{i-1}}Z$ and the set of $G^{0F}$-orbits on ${}^{\tt_i}Z$. We may assume 
that $i=1$. Let $\Ph$ be the fibre of ${}^1\vt$ at some $\FF_q$-rational point 
$(\tP,\tP',\ti\g)\in{}^{\tt_1}Z_{J_1,D}$. It is enough to show that $\Ph^F\ne\em$ and
that $\Ph^F$ is contained in a single $G^{0F}$-orbit. Since $\Ph$ is an affine space
(see \cite{\PCS, 3.12(b)}) defined over $\FF_q$, it must contain some $F$-fixed point
$(P,P',\g)$. By \cite{\PCS, 3.8}, $\Ph$ is a homogeneous $U_P\cap P'$ space (the action
being the restriction of the $G^0$-action on ${}^\tt Z$) and the isotropy group of
$(P,P',\g)$ is $U_P\cap U_{P'}$ (see \cite{\PCS, 3.9}). Since 
$U_P\cap P',U_P\cap U_{P'}$ are connected, it follows that $(U_P\cap P')^F$ acts 
transitively on $\Ph^F$; thus, $\Ph^F$ is contained in a single $G^{0F}$-orbit, as 
required.

Using this and the definitions, we see that (a) would be a consequence of the following
statement:

{\it Let $\ca'$ be a set of representatives for the isomorphism classes of admissible
simple perverse sheaves $A$ on $C$ such that $F^*A\cong A$. For each $A\in\ca'$ we 
choose an isomorphism $\a:F^*A@>\si>>A$; we define $A'$ (a simple perverse sheaf on 
$R$) as in 26.3 and let $\a':F^*A'@>\si>>A'$ be the isomorphism induced by $\a$. Then 
the characteristic functions $\c_{A',\a'}$ (one for each $A\in\ca'$) form a 
$\bbq$-basis of the vector space of functions $R^F@>>>\bbq$ that are constant on 
$G^{0F}$-orbits.}
\nl
This is an immediate consequence of 21.21 applied to $N_GL$ instead of $G$. The theorem
is proved.

\head 27. Induction\endhead
\subhead 27.1\endsubhead
Let $P$ be a parabolic of $G^0$ with Levi $L$. Let $G'=N_GP\cap N_GL$. Then $G'{}^0=L$.
Define a homomorphism $\x:N_GP@>>>G'$ by $\x(z\o)=z$ where $z\in G',\o\in U_P$ (see 
1.26). Let $D$ be a connected component of $G$ such that $N_GP\cap D\ne\em$. Then 
$N_GP\cap D$ is a connected component of $N_GP$ and $D'=G'\cap D$ is a connected 
component of $G'$. Consider the diagram
$$D'@<a<<V_1@>a'>>V_2@>a''>>D$$
where 

$V_1=\{(g,x)\in G\T G^0;x\i gx\in N_GP\cap D\}$,

$V_2=\{(g,xP)\in G\T G^0/P;x\i gx\in N_GP\cap D\}$,

$a(g,x)=\x(x\i gx),a'(g,x)=(g,xP),a''(g,xP)=g$.
\nl
Then $a,a'$ are smooth morphisms with connected fibres. To any $A\in\cm(D')$ which is 
$L$-equivariant for the conjugation action of $L$ on $D'$ we associate a complex 
$\ind_{D'}^D A\in\cd(D)$ as follows. The complex $\ta A\in\cm(V_1)$ is $P$-equivariant 
for the action $p:(g,x)\m(g,xp\i)$ of $P$ on $V_1$. Since $a'$ is a principal 
$P$-bundle, there is a well defined complex $A_1\in\cm(V_2)$ such that 
$\ta A=\ta' A_1$. We set $\ind_{D'}^D(A)=a''_!A_1$. 

\subhead 27.2\endsubhead
Let $L'$ be a Levi of a parabolic of $L$. Let $S$ be an isolated stratum of
$N_{G'}L'$ such that $S\sub D$ and $S$ normalizes some parabolic of $L$ with Levi $L'$.
Then $S$ is also an isolated stratum of $N_GL'$ such that $S$ normalizes some parabolic
of $G^0$ with Levi $L'$. We have a commutative diagram in which all squares 
except the top two are cartesian:
$$\CD
S'{}^*@<=<<S'{}^*@>=>>S'{}^*@.{}\\
@Ar_0AA @Ar_1AA @Ar_2AA@.\\
\hY'@<c<<\hZ_1@>c'>>\hZ_2@.{}\\
@Vq_0VV @Vq_1VV @Vq_2VV@.\\
\tY'@<b<<Z_1@>b'>>Z_2@<k<<\tY{}\\
@Vp_0VV @Vp_1VV @Vp_2VV@.\\
D'@<a<<V_1@>a'>>V_2@>a''>>D
\endCD$$
Here

$S'{}^*=\{g\in S;Z_{G'}(g_s)^0\sub L'\},
       \tY'=\{(g,lL')\in G'\T L/L';l\i gl\in S'{}^*\}$,

$S^*=\{g\in S;Z_G(g_s)^0\sub L'\},\tY=\{(g,xL')\in G\T G^0/L';x\i gx\in S^*\}$,

$Z_1=\{(g,lL',x)\in G\T\T L/L'\T G^0;l\i x\i gxl\in S'{}^*U_P\}$,

$Z_2=\{(g,xL'U_P)\in G\T G^0/(L'U_P); x\i gx\in S'{}^*U_P\}$,

$\hY'=\{(g,l)\in G'\T L;l\i gl\in S'{}^*\}$,

$\hZ_1=\{(g,l,x)\in G\T\T L\T G^0;l\i x\i gxl\in S'{}^*U_P\}$,

$\hZ_2=\{(g,xU_P)\in G\T G^0/U_P;x\i gx\in S'{}^*U_P\}$,

$c(g,l,x)=(\x(x\i gx),l),c'(g,l,x)=(g,xlU_P)$,

$b(g,lL',x)=(\x(x\i gx),lL'),b'(g,lL',x)=(g,xlL'U_P)$,

$k(g,xL')=(g,xL'U_P),q_0(g,l)=(g,lL'),q_1(g,l,x)=(g,lL',x)$,

$q_2(g,vU_P)=(g,vL'U_P),p_0(g,lL')=g,p_1(g,lL',x)=(g,x)$,

$p_2(g,vL'U_P)=(g,vP),r_0(g,l)=l\i gl$,
\nl
and $r_1(g,l,x)=s_1\in S'{}^*,r_2(g,vU_P)=s_2\in S'{}^*$ are defined by 

$l\i x\i gxl\in s_1U_P,v\i gv\in s_2U_P$.
\nl
Let $Y'=\cup_{l\in L}lS'{}^*l\i$ (a locally closed smooth irreducible subvariety of 
$D'$, see 3.16, 3.17). Let $Y'_1=a\i(Y')$; this is a locally closed smooth irreducible 
subvariety of $V_1$ which is $P$-stable since $Y'$ is stable under $L$-conjugacy (we 
use that $a$ is smooth with connected fibres). Hence $Y'_1=a'{}\i(Y'_2)$ where $Y'_2$ 
is a well defined locally closed smooth irreducible subvariety of $V_2$. Let
$Y=\cup_{x\in G^0}xS^*x\i$ (a locally closed smooth irreducible subvariety of $D$). Let
$\bY',\bY'_1,\bY'_2,\bY$ be the closure of $Y',Y'_1,Y'_2,Y$ in $D',V_1,V_2,D$. Then
$a'{}\i(\bY_2)=\bY_1=a\i(\bY')$. Let 

$\cw'=\{n\in N_LL';nSn\i=n\}/L'$, $\cw=\{n\in N_{G^0}L';nSn\i=n\}/L'$.
\nl
Then $\cw'\sub\cw$ are finite groups. Now $\cw'$ acts freely on $\tY',Z_1,Z_2$ by 

$n:(g,lL')\m(g,ln\i L'),n:(g,lL',x)\m(g,ln\i L',x)$,

$n:(g,xL'U_P)\m(g,xn\i L'U_P)$.
\nl
These actions are compatible with $\tY'@<b<<Z_1@>b'>>Z_2$. By 3.13, $p_0:\tY'@>>>Y'$ is
a principal $\cw'$-bundle. It follows that $p_1:Z_1@>>>Y'_1$ and $p_2:Z_2@>>>Y'_2$ are 
principal $\cw'$-bundles. 

Now $\cw$ acts freely on $\tY$ by $n:(g,xL')\m(g,xn\i L')$. From 3.13 we see that
$p:\tY@>>>Y,(g,xL')\m g$ is a principal $\cw$-bundle. From the definitions we see that 

(a) {\it the restriction of the $\cw$-action on $\tY$ to the subgroup $\cw'$ 
is compatible via $k:\tY@>>>Z_2$ with the $\cw'$-action on $Z_2$.}
\nl
Let $\ce\in\cs(S)$. Let $\tce$ be the local system on $\tY$ defined as in 5.6 (with 
$G,L',S$ instead of $G,L,S$). Let $\tce',\tce_1,\tce_2$ be the local system on 
$\tY',Z_1,Z_2$ (respectively) defined by 

$q_0^*\tce'=r_0^*\ce, q_1^*\tce_1=r_1^*\ce,q_2^*\tce_2=r_2^*\ce$.
\nl
We have $b^*\tce'=\tce_1=b'{}^*\tce_2$. Hence 
$(p_1)_!b^*\tce'=(p_1)_!b'{}^*\tce_2=(p_1)_!\tce_1$. Now $(p_0)_!\tce',(p_2)_!\tce_2$ 
are local systems on $Y',Y'_2$ and $a^*(p_0)_!\tce'=a'{}^*(p_2)_!\tce_2=(p_1)_!\tce_1$ 
as local systems on $Y'_1$. We have $\tce=k^*\tce_2$. Let

$K'=IC(\bY',(p_0)_!\tce')[\dim\tY'],K_1=IC(\bY'_1,(p_1)_!\tce_1)[\dim Z_1]$,

$K_2=IC(\bY'_2,(p_2)_!\tce_2)[\dim Z_2],K=IC(\bY,p_!\tce)[\dim\tY]$,
\nl
regarded as perverse sheaves on $D',V_1,V_2,D$, zero on 
$D'-\bY',V_1-\bY'_1,V_2-\bY'_2,D-\bY$ respectively. Since $a,a'$ are smooth morphisms 
with connected fibres, we see that $\ta K'=K_1=\ta'K_2$ in $\cm(V_1)$. Hence 

(b) $\ind_{D'}^D(K')=a''_!K_2$.
\nl
We show that

(c) $a''_!K_2=K$ {\it canonically}.
\nl
Let $P'$ be a parabolic of $L$ with Levi $L'$ such that $S\sub N_{G'}P'$. Then 
$P'_1=P'U_P$ is a parabolic of $G^0$ with Levi $L'$ such that $S\sub N_GP'_1$. We have 
a commutative diagram with cartesian squares
$$\CD
\tY'@<b<<Z_1@>b'>>Z_2\\
@Vj_0VV @Vj_1VV @Vj_2VV\\
X'@<e<<Z@>e'>>X\\
@Vt_0VV @Vt_1VV @Vt_2VV\\
D'@<a<<V_1@>a'>>V_2 
\endCD$$ 
Here

$X'=\{(g,lP')\in G'\T L/P';l\i gl\in\bS U_{P'}\}$,

$X=\{(g,vP'_1)\in G\T G^0/P'_1;v\i gv\in\bS U_{P'_1}\}$,

$Z=\{(g,lP',x)\in G\T\T L/P'\T G^0;l\i x\i gxl\in\bS U_{P'_1}\}$,

$e(g,lP',x)=(\x(x\i gx),lP'),e'(g,lP',x)=(g,xlP'_1)$,

$t_0(g,lP')=g,t_1(g,lP',x)=(g,x),t_2(g,vP'_1)=(g,vP)$,

$j_0(g,lL')=(g,lP'),j_1(g,lL',x)=(g,lP',x),j_2(g,xL'U_P)=(g,xP'_1)$. 
\nl 
Since $j_0$ is an open imbedding (by 5.5) we see that $j_1,j_2$ are open imbeddings. We
identify $\tY',Z_1,Z_2$ with open subsets of $X',Z,X$ via $j_0,j_1,j_2$. The 
composition $\tY@>k>>Z_2@>j_2>>X$ is the map $(g,xL')\m(g,xP'_1)$ which is an open 
imbedding by 5.5. Since $Z_2$ is an open subset of $X$ via $j_2$, we see that $\tY$ may
be identified with an open subset of $Z_2$ via $k$. Since $e,e'$ are smooth morphisms 
with connected fibres, we see that 

$\ti eIC(X',\tce')[\dim X']=IC(Z,\tce_1)[\dim Z]=\ti e'IC(X,\tce_2)[\dim X]$ in 
$\cm(Z)$.
\nl
From 5.7 we see that $K'=(t_0)_!IC(X',\tce')[\dim X']$. Hence

$K_1=(t_1)_!IC(Z,\tce_1)[\dim Z],K_2=(t_2)_!IC(X,\tce_2)[\dim X]$.
\nl
Since $\tce=\tce_2|_{\tY}$, we have $IC(X,\tce_2)=IC(X,\tce)$ and 
$K_2=(t_2)_!IC(X,\tce)[\dim X]$. The composition $a''t_2:X@>>>D$ is $(g,vP'_1)\m g$. 
Using 5.7 we have $a''_!K_2=(a''t_2)_!IC(X,\tce)[\dim X]=K$. This proves (c).

Define $\EE=\op_{w\in\cw}\EE_w$ as in 7.10 (for $G,L',S,\ce$ instead of $G,L,S,\ce$).
Define $\EE'=\op_{w\in\cw'}\EE'_w$ in the same way (for $G',L',S,\ce$ instead of 
$G,L,S,\ce$). Then $\EE$ is naturally an algebra and $\EE'$ is a subalgebra of $\EE$. 
Since $a,a'$ are smooth morphisms with connected fibres, we see that 

$\EE'=\End((p_0)_!\tce')=\End(K_0)=\End(K_1)=\End(K_2)$.
\nl
Now $a''_!$ defines a ring homomorphism $\End(K_2)@>>>\End(a''_!K_2)$. Thus $a''_!K_2$ 
becomes an $\EE'$-module. On the other hand $\EE=\End(K)$. Using (a) and the 
definitions we see that the restriction of the $\EE$-module structure of $K$ to $\EE'$ 
corresponds under (c) to the $\EE'$-module structure on $a''_!K_2$. 

Let $\G$ be a subgroup of $\cw'$ and let $\EE'_\G=\op_{w\in\G}\EE'_w$, a subalgebra of 
$\EE'$ hence of $\EE$. Let $\r$ be a $\EE'_\G$-module of finite dimension over $\bbq$. 
Let 

$K'(\r)=\Hom_{\EE'_\G}(\r,K')\in\cm(D'),K_1(\r)=\Hom_{\EE'_\G}(\r,K_1)\in\cm(V_1)$,

$K_2(\r)=\Hom_{\EE'_\G}(\r,K_2)\in\cm(V_2),K(\r)=\Hom_{\EE'_\G}(\r,K)\in\cm(D)$.
\nl
Then $\ta K'(\r)=K_1(\r)=\ta'K_2(\r)$ and
$$\ind_{D'}^D(K'(\r))=a''_!K_2(\r)=K(\r)\in\cm(D).\tag d$$

\subhead 27.3\endsubhead
Let $P,L,G',D,D'$ be as in 27.1. Let $Q$ be a parabolic of $L$ with Levi $M$ such that 
$N_{G'}Q\cap D'\ne\em$. Let $G''=N_{G'}Q\cap N_{G'}M$. Then $G''{}^0=M$,
$N_{G'}Q\cap D'$ is a connected component of $N_{G'}Q$ and $D''=G''\cap D'$ is a 
connected component of $G''$. Let $R=QU_P$, a parabolic of $G^0$ with Levi $M$. We have
$R\sub P$ and $N_GR\cap D\ne\em$. (Indeed, $D''\sub N_GR\cap D$; more precisely,
$N_GR\cap N_GM\cap D=D''$.) Note that $N_GR\cap N_GM$ contains $G''$ as a subgroup of 
finite index and both have $D''$ as a connected component.

Let $M'$ be a Levi of a parabolic of $M$. Let $S$ be an isolated stratum of
$N_{G''}M'$ such that $S\sub D$ and $S$ normalizes some parabolic of $M$ with Levi 
$M'$. Then $S$ is also an isolated stratum of $N_{G'}M'$ such that $S$ normalizes some 
parabolic of $L$ with Levi $M'$. Moreover, $S$ is an isolated stratum of $N_GM'$ such 
that $S$ normalizes some parabolic of $G^0$ with Levi $M'$. Let 

$S''{}^*=\{g\in S;Z_{G''}(g_s)^0\sub M'\}$,

$\tY''=\{(g,mM')\in G''\T M/M';m\i gn\in S''{}^*\}$,

$S'{}^*=\{g\in S;Z_{G'}(g_s)^0\sub M'\}$,

$\tY'=\{(g,lL')\in G'\T L/M';l\i gl\in S'{}^*\}$,

$S^*=\{g\in S;Z_G(g_s)^0\sub M'\}$,

$\tY=\{(g,xM')\in G\T G^0/M';x\i gx\in S^*\}$.
\nl
Let 

$Y''=\cup_{m\in M}mS''{}^*m\i,Y'=\cup_{l\in L}lS'{}^*l\i,Y=\cup_{x\in G^0}xS^*x\i$.
\nl
Define 

$p'':\tY''@>>>Y'',p':\tY'@>>>Y',p:\tY@>>>Y$
\nl
by the first projection. Let 

$\cw''=\{n\in N_MM';nSn\i=n\}/M'$, 

$\cw'=\{n\in N_LM';nSn\i=n\}/M'$, 

$\cw=\{n\in N_{G^0}M';nSn\i=n\}/M'$. 
\nl
Then $\cw''\sub\cw'\sub\cw$ are finite groups. As in 3.13, $p''$ is a principal 
$\cw''$-bundle, $p'$ is a principal $\cw'$-bundle, $p$ is a principal $\cw$-bundle. Let
$\ce\in\cs(S)$. Define a local system $\tce$ on $\tY$ as 5.6 (with $G,M',S$ instead of 
$G,L,S$). Define a local system $\tce'$ on $\tY'$ as 5.6 (with $G',M',S$ instead of 
$G,L,S$). Define a local system $\tce''$ on $\tY''$ as 5.6 (with $G'',M',S$ instead of 
$G,L,S$). Then $p''_!\tce'',p'_!\tce',p_!\tce$ is a local system on $Y'',Y',Y$
respectively. We regard 

$K''=IC(\bY'',p''_!\tce'')[\dim\tY'']$,

$K'=IC(\bY',p'_!\tce')[\dim\tY']$,

$K=IC(\bY,p_!\tce)[\dim\tY]$,
\nl
as perverse sheaves on $D'',D',D$ respectively, zero on $D''-\bY'',D'-\bY',D-\bY$.
Define $\EE=\op_{w\in\cw}\EE_w$ as in 7.10 (for $G,M',S,\ce$ instead of $G,L,S,\ce$).
Define in the same way $\EE'=\op_{w\in\cw'}\EE'_w$ (for $G',M',S,\ce$ instead of 
$G,L,S,\ce$) and $\EE''=\op_{w\in\cw''}\EE''_w$ (for $G'',M',S,\ce$ instead of 
$G,L,S,\ce$). Then $\EE$ is naturally an algebra, $\EE'$ is a subalgebra of $\EE$ and
$\EE''$ is a subalgebra of $\EE'$. We have naturally 

$\EE''=\End(K''),\EE'=\End(K'),\EE=\End(K)$.
\nl
Thus, $K'',K',K$ are naturally $\EE''$-modules. Let $\r$ be a finite dimensional 
$\EE''$-module over $\bbq$. Let 

$K''(\r)=\Hom_{\EE''}(\r,K'')\in\cm(D'')$,

$K'(\r)=\Hom_{\EE''}(\r,K')\in\cm(D')$,

$K(\r)=\Hom_{\EE''}(\r,K)\in\cm(D)$.
\nl
Applying 27.2(d) with $G',Q,M,M'$ instead of $G,P,L,L'$ and $\G=\cw''$ we have

$\ind_{D''}^{D'}(K''(\r))=K'(\r)\in\cm(D')$.
\nl
Applying 27.2(d) with $G,P,L,M'$ instead of $G,P,L,L'$ and $\G=\cw''$ we have

$\ind_{D'}^D(K'(\r))=K(\r)\in\cm(D)$.
\nl
Applying 27.2(d) with $G,R,M,M'$ instead of $G,P,L,L'$ and $\G=\cw''$ we have

$\ind_{D''}^D(K''(\r))=K(\r)\in\cm(D')$.
\nl
Hence we have the following {\it transitivity formula}:
$$\ind_{D'}^D(\ind_{D''}^{D'}(K''(\r)))=\ind_{D''}^D(K''(\r)).\tag a$$

\enddocument